\documentclass[11pt]{article}

\usepackage{amsmath}
\usepackage{amsmath, amsthm, amsfonts, amssymb, color}
\usepackage{mathrsfs}
\usepackage{latexsym,bm}
\usepackage{graphicx}

\newtheorem{theorem}{Theorem}[section]

\newtheorem{lemma}[theorem]{Lemma}
\newtheorem{corollary}[theorem]{Corollary}
\newtheorem{prop}[theorem]{Proposition}
\newtheorem{remark}[theorem]{Remark}

\def\FF{\mathcal F}

\def\eps{\varepsilon}

\definecolor{wco}{rgb}{0.5,0.2,0.3}

\numberwithin{equation}{section} 

\begin{document}

\allowdisplaybreaks

\title{\bf Law of large numbers for branching symmetric Hunt processes with
measure-valued branching rates }
\author{{\bf Zhen-Qing Chen\footnote{Department of Mathematics, University of Washington, Seattle,
WA 98195, USA. E-mail: zqchen@uw.edu}},
\quad {\bf Yan-Xia Ren}\footnote{LMAM School of Mathematical Sciences \& Center for
Statistical Science, Peking University, Beijing, 100871, P.R. China. E-mail: yxren@math.pku.edu.cn} \quad  \hbox{and} \quad {\bf Ting
Yang}
\footnote{Corresponding author. School of Mathematics and Statistics, Beijing Institute of Technology, Beijing, 102488, P.R.China. Email: yangt@bit.edu.cn}
\footnote{Beijing Key Laboratory on MCAACI, Beijing, 102488, P.R. China.}
}
\date{}
\maketitle

\begin{abstract}
We establish weak and strong law of large numbers for a class of branching symmetric Hunt processes with the branching rate being a smooth measure with respect to the underlying Hunt process, and the branching mechanism being  general and state-dependent. Our work is motivated by recent work
on strong law of large numbers for branching symmetric Markov processes by Chen-Shiozawa \cite{Chen&Shiozawa}
and for branching diffusions  by Engl\"ander-Harris-Kyprianou \cite{EHK}.
Our results can be applied to some interesting examples that are
covered by neither of these papers.
\end{abstract}

\medskip


\medskip


\section{Introduction}
\subsection{Motivation}
The law of large numbers (LLN)
has been the object of interest for measure-valued Markov processes including branching Markov processes and superprocesses.
For branching Markov processes, the earliest work in this field
dates back to 1970s when Watanabe \cite{Watanabe1,Watanabe2} studied
the asymptotic properties of a branching symmetric diffusion, using
a suitable Fourier analysis. Later, Asmussen and Hering
\cite{Asmussen} established an almost-sure limit theorem for a
general supercritical branching process under some regularity conditions.
Recently, there is a
revived interest in this field using modern techniques such as
Dirichlet form method, martingales and spine method.
Chen and Shiozawa \cite{Chen&Shiozawa} (see also
\cite{Shiozawa}) used a Dirichlet form and spectral theory approach
to obtain
strong law of large numbers (SLLN)
for branching symmetric Hunt processes. Among other assumptions, a spectral gap condition was used to obtain a Poincar\'{e} inequality that plays an important role in the proof of SLLN along lattice times.
They proved SLLN holds for branching processes under the assumptions that the branching rate is given by a
measure in Kato-class
$\mathbf{K}_{\infty}(X)$ and the branching mechanism has bounded second moment.

The spine method developed recently for measure-valued Markov processes is a powerful probabilistic tool in studying various properties of these processes;
see, e.g., \cite{EKW,E,EHK,EK,HH,LRSa,LRSb}.
In \cite{EHK},  Engl\"{a}nder,
Harris and Kyprianou used spine method to obtain
SLLN for branching (possibly non-symmetric) diffusions
corresponding to the operator $Lu+\beta(u^2-u)$ on
a domain $D\subset\mathbb{R}^d$
(where $\beta\ge 0$ is non-trivial) under certain spectral conditions.
They imposed a condition on how far  particles can spread in space
(see condition (iii) on page 282 of \cite{EHK}). That the underlying process is a diffusion plays an important role in their argument
and the branching rate there has to be a function rather than a measure.
The approach of \cite{EHK} also involves $p$-th moment  calculation with
$p>1$ which may not be valid for general branching mechanisms.
Recently, Eckhoff,  Kyprianou and Winkel \cite{EKW}  discussed the strong law of large numbers (SLLN)  along lattice times for branching diffusions, which serves as the backbone or skeleton for superdiffusions.
It is proved in \cite[Theorem 2.14]{EKW} that,
if the branching mechanism satisfies a $p$-th moment condition with $p\in (1,2]$, the underlying diffusion and  the
support of the branching diffusion satisfy conditions similar to that presented in  \cite{EHK},  then SLLN along lattice times holds.

In this paper, we combine
the functional analytic methods used in \cite{Chen&Shiozawa} with spine techniques
to study weak and strong laws of
large numbers for
branching symmetric Hunt processes as well as the $L\log L$ criteria.
This approach allows us to obtain new results for a large class of
branching Markov processes,
for which
(i) the underlying spatial motions are general symmetric Hunt processes, which can be discontinuous
and may not be
intrinsically ultracontractive;
(ii) the branching rates are given by general smooth measures rather than functions or Kato class measures;
(iii) the offspring distributions are only assumed to have bounded
first moments with no assumption on their second moments.
In addition, we use  $L^{1}$-approach instead of $L^{p}$-approach
for $p\in(1,2]$.
Now we describe the setting and main results of this paper in detail, followed by several examples illuminating the main results.

\subsection{Branching symmetric Hunt processes and assumptions}
Suppose we are given three initial ingredients: a Hunt process, a smooth measure and a branching mechanism.
We introduce them one by one:
\begin{itemize}
 \item {{\it A Hunt process $X$:}  Suppose $E$ is a locally compact separable metric space and
 $E_{\partial} := E \cup \{ \partial \}$  is
 its one point compactification.
 $m$ is a positive Radon measure on $E$
 with full support. Let $X=(\Omega,\mathcal{H},\mathcal{H}_{t},\theta_{t},X_{t},\Pi_{x},\zeta)$ be a $m$-symmetric Hunt process on $E$.
 Here $\{\mathcal{H}_{t}:\ t\ge 0\}$ is the minimal admissible filtration, $\{\theta_{t}:\ t\ge
 0\}$ the time-shift operator of $X$ satisfying
 $X_{t}\circ\theta_{s}=X_{t+s}$ for $s,t\ge 0$, and $\zeta:=\inf\{t>0:\
 X_{t}=\partial\}$ the life time of $X$.
  Suppose for each $t>0$, $X$ has a
  symmetric transition density function
  $p(t,x,y)$ with respect to the measure $m$.
  Let $\{P_{t}:\ t\ge 0\}$ be
 the Markovian transition semigroup of $X$, i.e.,
 $$P_{t}f(x):=\Pi_{x}\left[f(X_{t})\right]=\int_{E}f(y)p(t,x,y)m(dy)$$
 for any non-negative measurable function $f$. The symmetric
 Dirichlet form on $L^{2}(E,m)$ generated by $X$ will be denoted as
 $(\mathcal{E},\mathcal{F})$:
 $$
 \mathcal{F}= \Big \{u\in L^{2}(E,m):\ \lim_{t\to 0}\frac{1}{t}\int_{E}\left(u(x)-P_{t}u(x)\right)u(x)m(dx)<+\infty \Big\},
 $$
 $$\mathcal{E}(u,v)=\lim_{t\to 0}\frac{1}{t}\int_{E}\left(u(x)-P_{t}u(x)\right)v(x)m(dx),\quad u,v\in\mathcal{F}.$$
 It is known (cf. \cite{CF})
 that $(\mathcal{E},\mathcal{F})$ is quasi-regular and hence is quasi-homeomorphic to a regular Dirichlet form on a locally compact separable
 metric space. In the sequel, we assume that $X$ is $m$-irreducible in the sense that if $A\in\mathcal{B}(E)$ has positive $m$-measure, then $\Pi_{x}(T_{A}<+\infty)>0$ for
 all $x\in E$, where $T_{A}:=\inf\{t>0:\ X_{t}\in A\}$ is the first hitting time of $A$. }
 \item{{\it A branching rate $\mu$:} Suppose $\mu$ is a positive smooth Radon measure on $(E,\mathcal{B}(E))$.
 It
 uniquely determines a positive
 continuous additive functional (PCAF) $A^{\mu}_{t}$ by the following Revuz formula:
 $$\int_{E}f(x)\mu(dx)=\lim_{t\to 0}\frac{1}{t}\Pi_{m}\left[ \int_{0}^{t}f(X_{s})dA^{\mu}_{s}\right],
 \quad f\in \mathcal{B}^{+}(E).$$
 Here $\Pi_{m}(\cdot):=\int_{E}\Pi_{x}(\cdot)m(dx)$.}
 \item{{\it Offspring distributions
  $\{\{p_n(x): n\ge 0\},
  x\in E\}$:}
 Suppose
  $\{\{p_n(x): n\ge 0\},
  x\in E\}$
 is a family of probability mass functions such that $0\le p_n(x)\le 1$ and $\sum^\infty_{n=0}p_n(x)=1$.
  For each $x\in E$,
    $\{p_n(x): n\ge 0\}$
 serves as the offspring distribution of a particle  located at $x$.
 Let $\{A(x):x\in E\}$ be a collection of random variables taking values in $\{0,1,2,\cdots\}$ and distributed as $\mathrm{P}(A(x)=n)=p_{n}(x)$.
    }
\end{itemize}

Define
\begin{equation}
Q(x):=\sum_{n=0}^{+\infty}np_{n}(x),\quad x\in E.
\end{equation}
Throughout this paper we assume that the offspring distribution
$\{p_n(x): n\ge 0\}$
satisfies the following condition:
\begin{equation}
p_{0}(x)\equiv 0,\quad p_{1}(x)\not\equiv 1 \quad
\hbox{and} \quad \sup_{x\in E}Q(x)<\infty.
\end{equation}

From these ingredients we can build a branching Markov process according to the following recipe: under a probability
measure $\mathbb{P}_{x}$, a particle starts from $x\in E$ and moves
around in $E_{\partial}$ like a copy of $X$. We use $\emptyset$ to
denote the original particle, $X_{\emptyset}(t)$ its position at
time $t$ and $\zeta_{\emptyset}$ its
fission time.
We say that
$\emptyset$ splits at the rate $\mu$ in the sense that
$$
\mathbb{P}_{x}\left(\zeta_{\emptyset}>t|X_{\emptyset}(s):s\ge 0\right)=\exp (-A^{\mu}_{t} ).
$$
When $\zeta_{\emptyset}\ge \zeta$, it dies at time $\zeta$. On the
other hand, when $\zeta_{\emptyset}<\zeta$, it splits into a random
number of children, the number being distributed as a copy of
$A(X_{\emptyset}(\zeta_{\emptyset}-))$. These children, starting
from their point of creation, will move and reproduce
independently in the same
way as their parents. If a particle $u$ is alive at time $t$, we refer to its location in $E$ as $X_{u}(t)$. Therefore the time-$t$ configuration is a $E$-valued point process $\mathbb{X}_{t}=\{X_{u}(t): u\in
 {\cal Z}_t \}$,
 where ${\cal Z}_t$  is the set of particles alive at time $t$. With abuse of  notation, we can also regard $\mathbb{X}_{t}$ as  a random point measure on $E$ defined by
$\mathbb{X}_{t}:=\sum_{u\in{\cal Z}_t}\delta_{X_{u}(t)}$.
Let  $({\cal F}_t)_{t\ge 0}$ be the natural filtration of $\mathbb{X}$ and ${\cal F}_\infty=
\sigma\{\mathcal{F}_{t}: t\ge 0\}$.
Hence it defines a branching symmetric Hunt process
$\mathbb{X}=(\mathbf{\Omega},\mathcal{F}_{\infty},\mathcal{F}_{t},\mathbb{X}_{t},\mathbb{P}_{x})$ on $E$ with the motion component $X$,  the branching rate measure $\mu$ and the branching mechanism function $\{p_{n}(x):n\ge 0\}$.
When the branching rate measure $\mu$ is absolutely continuous with respect to $m$, i.e., $\mu(dy)=\beta(y)m(dy)$ for some non-negative function $\beta$, the corresponding PCAF $A^{\mu}_{t}$ is equal to $\int_{0}^{t}\beta(X_{s})ds$, and given that a particle $u$ is alive at $t$, its probability of splitting in $(t,t+dt)$ is $\beta(X_{u}(t))dt+o(dt)$. Since the function $\beta$ determines the rate at which every particle splits,
$\beta$ is  called the branching rate function in literature.

Throughout this paper we use $\mathcal{B}_{b}(E)$ (respectively,
$\mathcal{B}^{+}(E)$) to denote
the space of bounded (respectively,
non-negative) measurable functions on $(E,\mathcal{B}(E))$.
Any function $f$ on $E$
will be automatically extended to $E_{\partial}$ by
setting $f(\partial)=0$. We use $\langle f,g\rangle$ to denote
$\int_{E}f(x)g(x)m(dx)$ and
``$:=$"
as a way of definition.
For $a,b\in \mathbb{R}$,
$a\wedge b:=\min\{a,b\}$,
$a\vee b:=\max\{a,b\}$,
 and $\log^{+}a:=\log (a\vee 1)$.

For every $f\in\mathcal{B}^+(E)$ and $t\ge 0$, define
$$\mathbb{X}_{t}(f):=\sum_{u\in{\cal Z}_{t}}f(X_{u}(t)).$$
 We define the Feynman-Kac semigroup $P^{(Q-1)\mu}_{t}$ by
 \begin{equation}
 P^{(Q-1)\mu}_{t}f(x):=\Pi_{x}\left[
 \exp(A^{(Q-1)\mu}_{t})
 f(X_{t})\right],\quad
 f\in\mathcal{B}^{+}(E).\label{FK}
 \end{equation}
 Since $X$ has a transition density function,
 it follows
 that for each $t>0$, $P^{(Q-1)\mu}_{t}$ admits
 an integral kernel
 with respect to the measure $m$. We denote this kernel by
 $p^{(Q-1)\mu}(t,x,y)$.
  The semigroup $P^{(Q-1)\mu}_{t}$ associates with
 a quadratic form $(\mathcal{E}^{(Q-1)\mu},\mathcal{F}^{\mu})$,  where
$ \mathcal{F}^{\mu}  =  \mathcal{F} \cap L^2(E ; {\mu} )$ and
$$
\mathcal{E}^{(Q-1)\mu}(u,u)  =  \mathcal{E}(u,u)-\int_{E}u(x)^{2}(Q(x)-1)\mu(dx),
 \quad   u\in \mathcal{F}^{\mu}.
$$
We say that a signed smooth measure $\nu$
belongs to the Kato class $\mathbf{K}(X)$,
if \begin{equation}
\lim_{t\downarrow 0}\sup_{x\in
E}\Pi_{x}[A^{ | \nu |}_{t}]=0.\label{Kato}
\end{equation}
For every non-negative $\nu\in\mathbf{K}(X)$,
we have $\|G_\alpha \nu \|_\infty < \infty $ for  every $\alpha>0$,
where $G_\alpha$ is the $\alpha$-resolvent of $X$ and
 $G_\alpha \nu $ is the $\alpha$-potential of $\nu $.
Define  $\mathcal{E}_{\alpha}(u,u):=\mathcal{E}(u,u)+\alpha \int_{E}u^{2}dm$.
By \cite{SV},
\begin{equation}  \label{0.1}
 \int_{E}u(x)^{2}\nu(dx)\le
 \| G_\alpha \nu \|_\infty \mathcal{E}_{\alpha}(u,u),
 \quad  u\in\mathcal{F}.
\end{equation}
Thus when $\mu$ is in $\mathbf{K}(X)$,
 $\mathcal{F}^{\mu}=\mathcal{F}$, and
the quadratic form $(\mathcal{E}^{(Q-1)\mu},\mathcal{F}^{\mu})$ is
bounded from below.

For an arbitrary smooth measure $\mu$,
we define
\begin{equation}
\lambda_{1}:=\inf\left\{\mathcal{E}^{(Q-1)\mu}(u,u):u\in
\mathcal{F}^{\mu} \hbox{ with }  \int_{E}u(x)^{2}m(dx)=1\right\}.\label{l1}
\end{equation}

\noindent\textbf{Assumption 1.}
Let $\mu$ be a non-negative smooth measure on $E$ so that
the Schr\"odinger semigroup $P^{(Q-1)\mu}_t$ admits a symmetric kernel
$p^{(Q-1)\mu}(t, x, y)$ with respect to the measure $m$ and is jointly continuous
in $(x, y)\in E\times E$ for every $t>0$. Moreover,
$\lambda_{1}\in (-\infty,0)$ and
there is a positive continuous function $h \in \mathcal{F}^{\mu}$ with
$\int_{E}h(x)^{2}m(dx)=1$
so that
$ \mathcal{E}^{(Q-1)\mu}(h, h) = \lambda_1$.

\medskip

Observe that if $u$ is a minimizer for \eqref{l1}, then so is $| u|$.
Assumption 1 says that there is a minimizer  for \eqref{l1} that can be chosen
to be positive
and continuous. Clearly the following property holds for $h$:
\begin{equation}
\mathcal{E}(h,v)=\int_{E}h(x)v(x)(Q(x)-1)\mu(dx)+\lambda_{1}\langle
h,v\rangle  \quad \hbox{for every }  v\in\mathcal{F}^{\mu}.\label{1.3}
\end{equation}

The finiteness of $\lambda_{1}$ implies that
the bilinear form $(\mathcal{E}^{(Q-1)\mu},\mathcal{F}^{\mu})$ is
bounded from below, and hence by \cite{A&M},
$\{P^{(Q-1)\mu}_{t}:t\ge 0\}$ is a strongly continuous semigroup on
$L^{2}(E,m)$. Let
$\sigma(\mathcal{E}^{(Q-1)\mu})$ denote the spectrum of the self-adjoint operator associated with $\mathcal{E}^{(Q-1)\mu}$.
Let $\lambda_{2}$ be the second bottom of $\sigma(\mathcal{E}^{(Q-1)\mu})$, that is,
$$
\lambda_{2}:=\inf \Big\{ \mathcal{E}^{(Q-1)\mu}(u,u):\ u\in\mathcal{F}^{\mu}, \ \int_{E}u(x)h(x)m(dx)=0,\ \int_{E}u(x)^{2}m(dx)=1
\Big\}.
$$

\medskip

\noindent\textbf{Assumption 2.} There is a positive spectral gap in
$\sigma(\mathcal{E}^{(Q-1)\mu})$:
$\lambda_{h}:=\lambda_{2}-\lambda_{1} >0$.
\medskip

Define $h$-transformed semigroup $\{P^h_t; t\geq 0\}$ from
$\{ P^{(Q-1)\mu}_{t}; t\geq 0\}$ by
\begin{equation}\label{e:1.8}
P^{h}_{t}f(x)=
\frac{e^{\lambda_{1}t}}{h(x)}P^{(Q-1)\mu}_{t}(hf)(x)
\quad\mbox{ for }x\in E\mbox{ and }f\in\mathcal{B}^{+}(E).
\end{equation}
Then it is easy to see that $\{P^h_t; t\geq 0\}$ is an $\widetilde m$-symmetric semigroup, where $\widetilde m:=h^2 m$, and
$1$ is an eigenfunction of $P^{h}_{t}$
with eigenvalue $1$. Furthermore the spectrum of the infinitesimal generator of
$\{P^h_t; t\geq 0\}$ in $L^2(E; \widetilde m)$ is the spectrum of the infinitesimal generator of $\{ P^{(Q-1)\mu}_{t}; t\geq 0\}$ in $L^2(E; m)$ shifted by
$\lambda_1$. Hence under Assumption 2,
we have the following Poincar\'{e} inequality:
\begin{equation} \label{p1}
\|P^{h}_{t}\varphi\|_{L^{2}(E,\widetilde{m})}\le
e^{-\lambda_h t}\|\varphi\|_{L^{2}(E,\widetilde{m})}
\end{equation}
for all $\varphi\in L^{2}(E,\widetilde{m})$ with
$\int_{E}\varphi(x)\widetilde{m}(dx)=0$.

\begin{remark}\rm
If the underlying process $X$ satisfies that for each $t>0$, the transition density function $p(t,x,y)$ is bounded and is continuous in $x$ for every fixed $y\in E$
and that the branching rate measure $\mu$ is in
the Kato class $\mathbf{K}(X)$ of $X$,
then it follows from \cite{ABM} that the Feymann-Kac semigroup $P^{(Q-1)\mu}_{t}$ maps bounded functions to continuous functions and is bounded from $L^{p}(E,m)$ to $L^{q}(E,m)$ for any $1\le p\le q\le +\infty$. By Friedrichs theorem,
Assumptions 1 and 2 hold if we assume in addition that
\begin{center}
the embedding of $(\mathcal{F},\mathcal{E}_{1})$ into $L^{2}(E,\mu)$ is compact.
\end{center}
Such an assumption is imposed in \cite{Chen&Shiozawa} to
ensure the spectral gap condition and to obtain Poincar\'{e} inequality \eqref{p1}.
\end{remark}

\subsection{Main results}

Recall that  $({\cal F}_t)_{t\ge 0}$ is the natural filtration of  $\mathbb{X}$.
Observe that (cf. \cite[Lemma 3.3]{Shiozawa})
for every $x\in E$ and every
$f\in\mathcal{B}^{+}(E)$,
\begin{equation}
\mathbb{E}_{x}\left[\mathbb{X}_{t}(f)\right]=P^{(Q-1)\mu}_{t}f(x).\label{manytoone}
\end{equation}
 It is easy to see that
 $M_{t}:=e^{\lambda_{1}t}\mathbb{X}_{t}(h)$ is a positive $\mathbb{P}_{x}$-martingale with respect to $\mathcal{F}_{t}$. Let $M_{\infty}:=\lim_{t\to+\infty}M_{t}$.
It is natural to ask
when $M_{\infty}$ is non-degenerate, that is, when $\mathbb{P}_{x}(M_{\infty}>0)>0$ for $x\in E$?
Under the assumptions that (i) $m(E)<\infty$,
(ii) the Feymann-Kac semigroup $P^{(Q-1)\mu}_{t}$ is intrinsically ultracontractive,
and (iii) $h$ is bounded,
it is proved in \cite{LRSb} that the condition
\begin{equation}\label{llogl'}
\int_{E}\int_{E}\sum_{k=0}^{+\infty}kp_{k}(y)
h(y)^{2} \log^{+}(k h(y))\mu(dy)<+\infty
\end{equation}
 is  necessary and sufficient for $M_{\infty}$
 to be non-degenerate.
Condition \eqref{llogl'} is usually called the
$L\log L$ criteria.
The first main result of this paper reveals that,
in general, condition \eqref{llogl} below is sufficient for $M_{\infty}$
  to be non-degenerate.

\begin{theorem}\label{them1}
 Suppose Assumptions 1-2 hold. If
\begin{equation}\label{llogl}
\int_{E} h(y)^{2} \log^{+}h(y)m(dy)+\int_{E}\sum_{k=0}^{+\infty}kp_{k}(y)
h(y)^{2} \log^{+}(k h(y))\mu(dy)<+\infty,
\end{equation}
   then $M_{t}$ converges to $M_{\infty}$ in $L^{1}(\mathbb{P}_{x})$ for every $x\in E$, and, consequently, $\mathbb{P}_{x}(M_{\infty}>0)>0$.
\end{theorem}

Thus
under condition \eqref{llogl}, $\mathbb{X}_{t}(h)$ grows  exponentially with rate $-\lambda_1$.
Note that when $h$ is bounded, \eqref{llogl} is equivalent to \eqref{llogl'}.
The next question to ask is that, for a general
test function
$f\in\mathcal{B}^{+}(E)$, what is the limiting behavior
of $ \mathbb{X}_{t}(f)$ as $t\to\infty$? By \eqref{e:1.8} and \eqref{manytoone}, it is not hard to deduce (see \eqref{p3} below) that for every $f\in\mathcal{B}^{+}(E)$
with $f\leq c h$ for some constant $c>0$,
 $$e^{\lambda_{1}t}\mathbb{E}_{x}\left[\mathbb{X}_{t}(f)\right]=h(x)P^{h}_{t}(f/h)(x)
 \to h(x)\langle f,h\rangle\quad\mbox{ as }t\to+\infty.$$
 So, the mean of $ \mathbb{X}_{t}(f)$ also grows  exponentially with rate $-\lambda_1$. Our previous question is related to the question: for $f\in\mathcal{B}^{+}(E)$ with
 $f\leq c h$ for some constant $c>0$,
 does $\mathbb{X}_{t}(f)$ grow
exponentially with the same rate?
If so, can one identify its limit?
We first answer these questions in Theorem \ref{them2} and
Corollary \ref{C:1.4}
in terms of  convergence  in $L^{1}(\mathbb{P}_{x})$ and in probability,
under an additional condition \eqref{W}.

Note that under Assumption 1, for every $t>0$, $P^{h}_{t}$ has a
symmetric
continuous transition density function $p^{h}(t,x,y)$ on $E\times E$
with respect to the measure $\widetilde{m}$, which is related to $p^{(Q-1)\mu}(t,x,y)$ by the following formula:
$$p^{h}(t,x,y)=e^{\lambda_{1}t}\frac{p^{(Q-1)\mu}(t,x,y)}{h(x)h(y)},
\quad  x,y\in E.$$

\begin{theorem}[Weak law of large numbers]\label{them2}
Suppose Assumptions 1-2 and \eqref{llogl} hold.
If there exists some $t_{0}>0$ such that
\begin{equation}\label{W}
\int_{E} p^{(Q-1)\mu} (t_{0}, y, y)m(dy)<+\infty,\quad \mbox{or equivalently,}
\quad\int_{E} p^h (t_{0}, y, y)\widetilde{m}(dy)<+\infty,
\end{equation}
then for any $x\in E$ and any $f\in\mathcal{B}^{+}(E)$ with $f\leq ch$ for some $c>0$,
we have
$$\lim_{t\to+\infty}e^{\lambda_{1}t}\mathbb{X}_{t}(f)= M_{\infty}\langle f,h\rangle\quad\mbox{  in }L^{1}(\mathbb{P}_{x}).$$
\end{theorem}

\begin{corollary}\label{C:1.4}
Under the assumptions of Theorem \ref{them2},
it holds that
$$\lim_{t\to+\infty}\frac{\mathbb{X}_{t}(f)}{\mathbb{E}_{x}[\mathbb{X}_{t}(f)]}=\frac{M_{\infty}}{h(x)}\quad\mbox{ in probability with respect to }\mathbb{P}_{x},$$
for every $x\in E$ and every $f\in\mathcal{B}^{+}(E)$ with
$f\leq ch$ for some $c>0$.
\end{corollary}

\medskip

For almost sure convergence result,
we need a stronger condition \eqref{0.2} below.

\begin{theorem}[Strong law of large numbers]\label{them3}
Suppose Assumptions 1-2 and \eqref{llogl} hold.
If there exists
$t_{1}>0$
such that
\begin{equation}\label{0.2}
\sup_{y\in E}\frac{p^{(Q-1)\mu} (t_{1}, y, y)}{h(y)^{2}}<+\infty,\quad\mbox{ or, equivalently,}
\quad\sup_{y\in E} p^{h} (t_{1}, y, y)<+\infty,
\end{equation}
then there exists $\mathbf{\Omega_{0}}\subset \mathbf{\Omega}$ of
$\mathbb{P}_{x}$-full probability for
every $x\in E$,
such that, for
every $\omega\in\mathbf{\Omega_{0}}$ and every $f\in
\mathcal{B}_{b}(E)$ with compact support whose set of discontinuous
points has zero $m$-measure, we have
\begin{equation}
\lim_{t\to+\infty}e^{\lambda_{1}t}\mathbb{X}_{t}(f)(\omega)=M_{\infty}(\omega)\langle
f,h\rangle.\label{continuous}
\end{equation}
\end{theorem}

\medskip

\begin{corollary}
Suppose the assumptions of
Theorem \ref{them3} hold and let $\mathbf{\Omega_{0}}$ be
defined in Theorem \ref{them3}.
Then
$$\lim_{t\to+\infty}\frac{\mathbb{X}_{t}(f)(\omega)}{\mathbb{E}_{x}\left[\mathbb{X}_{t}(f)\right]}=\frac{M_{\infty}(\omega)}{h(x)}$$
for every $\omega\in\mathbf{\Omega_{0}}$ and for every $f\in
\mathcal{B}_{b}(E)$ with compact support whose set of discontinuous
points has zero $m$-measure.
\end{corollary}

\begin{remark}\rm
The condition \eqref{W}
is equivalent to
$$\int_{E}\int_{E}p^{h}(t_{0}/2,x,y)^{2}\widetilde{m}(dy)\widetilde{m}(dx)<+\infty.$$
Hence by \cite[Page 156]{Schwartz}, $P^{h}_{t}$ is a compact operator on $L^{2}(E,\widetilde{m})$ for every $t\geq t_0/2$.
 Consequently Assumption 2 is automatically satisfied if either \eqref{W} or \eqref{0.2} holds.
\end{remark}

To understand condition \eqref{0.2}, we  give some equivalent statements of \eqref{0.2} under our Assumptions 1-2.

\begin{prop}\label{P:1.8} Suppose Assumptions 1-2 hold.
The following are equivalent to \eqref{0.2}.
\begin{description}
\item{\rm (i)} There exists $t_1>0$ such that for any $t>t_{1}$,
\begin{equation}
\sup_{x,y\in E} |
p^{h}(t,x,y)
-1|
\le c_{1}e^{-c_{2}t}\label{1.17}
\end{equation}
for some $c_{1},c_{2}>0$.

\item{\rm (ii)}
 \begin{equation}\label{0.3}
p^{h}(t,x,y)
\to 1,\mbox{ as }t\uparrow +\infty\mbox{ uniformly in }(x,y)\in E\times
E.
\end{equation}

\item{\rm (iii)} There exist constants $t, c_{t}>0$ such that
\begin{equation}
\label{AIUC'}
p^{(Q-1)\mu}(t,x,y)
\le c_{t}h(x)h(y)\quad\hbox{for every } x,y\in
E.
\end{equation}
\end{description}
\end{prop}

Property \eqref{AIUC'}
is called \textit{asymptotically intrinsically ultracontractive} (AIU)
by Kaleta and L\H{o}rinczi in \cite{Kaleta}.
If
the inequality in \eqref{AIUC'}
is true for every $t>0$, and every $x,y\in E$,
then
$\{P^{(Q-1)\mu}_{t}: t>0\}$
is
called \textit{intrinsically ultracontractive}
(IU).
 It is shown in
\cite{Kaleta} that in case of symmetric $\alpha$-stable processes
($\alpha\in (0,2)$),
(AIU) is a weaker property than (IU).

\subsection{Examples}
In this subsection,
we illustrate our main results by several concrete examples.
For simplicity,
we consider binary branching only, i.e., every particle gives birth to
precisely two children, in which case $Q(x)\equiv 2$ on $E$.
Since $\sum_{k=0}^{+\infty}kp_{k}(x)\log^{+}k\equiv 2\log 2$ on $E$, condition \eqref{llogl} is reduced to
\begin{equation}
\int_{E} h(y)^{2}  \log^{+}h(y)\left(m(dy)+\mu(dy)\right)<+\infty.\label{llogl1}
\end{equation}

\noindent\textbf{Example 1. [WLLN for branching OU processes with a quadratic
branching rate function]}  Let $E=\mathbb{R}^{d}$. In Example 10 of
\cite{EHK}, $(X,\Pi_{x})$ is an
Ornstein-Ulenbeck (OU) process
on $\mathbb{R}^{d}$ with
infinitesimal generator
$$\mathcal{L}=\frac{1}{2}\sigma^{2}\Delta-cx\cdot\nabla\quad\mbox{  on }\mathbb{R}^{d},$$
where $\sigma,\ c>0$. Without loss of generality, we assume
$\sigma=1$.
Let $m(dx)=\left(\frac{c}{\pi}\right)^{d/2}e^{-c|x|^{2}}dx$. Then $X$ is symmetric with respect
to the probability measure $m$, and
the Dirichlet form $(\mathcal{E},\mathcal{F})$ of $X$ on
$L^{2}(\mathbb{R}^{d},m)$ is given by
$$\mathcal{F}=\left\{f\in L^{2}(\mathbb{R}^{d},m):\ \int_{\mathbb{R}^{d}}|\nabla f(x)|^{2}m(dx)<+\infty\right\},$$
$$\mathcal{E}(u,u)=\frac{1}{2}\int_{\mathbb{R}^{d}}|\nabla f(x)|^{2}m(dx).$$
Let $\beta(x)=b|x|^{2}+a$ with $a,b>0$ be the branching rate function. Let $P^{\beta}_{t}$ be the corresponding Feynman-Kac semigroup,
$$P^{\beta}_{t}f(x):=\Pi_{x}\left[
\exp\left(\int_{0}^{t}\beta(X_{s})ds\right)
f(X_{t})\right].$$
Suppose $c>\sqrt{2b}$ and $\alpha=\sqrt{c^{2}-2b}$. Let
$$\lambda_{c}:=\inf\{\lambda\in\mathbb{R}:\ \mbox{there exists }u>0 \mbox{ such that }(\mathcal{L}+\beta-\lambda)u=0\mbox{ in }\mathbb{R}^{d}\}$$
be the generalized principal eigenvalue. Let $\phi$ denote the corresponding ground state, i.e., $\phi>0$ such that
$(\mathcal{L}+\beta-\lambda_{c})\phi=0$.
As  indicated
in \cite{EHK}, $\lambda_{c}=\frac{1}{2}(c-\alpha)+a>0$ and
$\phi(x)=\left(\frac{\alpha}{c}\right)^{d/4}\exp(\frac{1}{2}(c-\alpha)|x|^{2})$. Note that $\phi\in\mathcal{F}^{\mu}$ and $\phi=e^{-\lambda_{c}t}P^{\beta}_{t}\phi$ on $\mathbb{R}^{d}$. It is easy to see that in this example $\lambda_{1}=-\lambda_{c}$ and $h(x)=\phi(x)$.
The transformed process $(X^{h}, \Pi^{h}_{x})$ is also an OU process with
infinitesimal generator
$$\mathcal{L}=\frac{1}{2}\Delta-\alpha x\cdot\nabla\mbox{  on }\mathbb{R}^{d}.$$
Note that its invariant probability measure is
$\widetilde{m}(dx)=h(x)^{2}m(dx)=\left(\frac{\alpha}{\pi}\right)^{d/2}e^{-\alpha|x|^{2}}dx$.
Let $p^{h}(t,x,y)$ be the transition density
of $X^{h}$ with respect to $\widetilde{m}$. It is known that
\begin{eqnarray}p^{h}(t,x,y)&=&\left(\frac{1}{1-e^{-2\alpha t}}\right)^{d/2}\exp\left( -\frac{\alpha}{(e^{2\alpha t}-1)}\left(|x|^{2}+|y|^{2}-2x\cdot y e^{\alpha t}\right)\right) \nonumber
\end{eqnarray}
In particular,
$$
p^{h}(t,x,x)=\left(\frac{1}{1-e^{-2\alpha t}}\right)^{d/2}\exp\left( \frac{2\alpha}{e^{\alpha t}+1}|x|^{2}\right).
$$
Thus
$\int_{\mathbb{R}^{d}}p^{h}(t,x,x)\widetilde{m}(dx)<+\infty$ for $t>0$.
Moreover, we observe that
condition \eqref{llogl1} is satisfied for this example.
Therefore Theorem \ref{them2} holds for this example.

This example does not satisfy the assumptions in \cite{Chen&Shiozawa}. To be more specific,
here the ground state $h$ is unbounded and $\beta (x)= b |x|^2 + a$ is not in
the Kato class ${\bf K}_\infty (X)$ of $X$.

\bigskip
\noindent\textbf{Example 2. [WLLN for branching Hunt processes with
a bounded branching rate function]} Let $E$ be a locally compact separable metric space and $m$ a positive Radon measure on $E$ with full support. Suppose the branching rate function
$\beta$ is a non-negative bounded function on $E$. Suppose the underlying
Hunt process $(X,\Pi_{x})$ satisfies that for every $t>0$, there exists a family of
jointly continuous, symmetric and positive kernels $p(t,x,y)$ such
that $P_{t}f(x)=\int_{E}p(t,x,y)f(y)m(dy)$, and
that there exists $s_{1}>0$ so that
\begin{equation}
\int_{E}p(s_{1},x,x)m(dx)<+\infty.\label{eg1.1}
\end{equation}
In this case the Feyman-Kac semigroup
$$P^{\beta}_{t}f(x):=\Pi_{x}\left[
\exp\left(\int_{0}^{t}\beta(X_{s})ds\right)
f(X_{t})\right]$$
has a jointly continuous and positive kernel $p^{\beta}(t,x,y)$.
It is easy to see that
\begin{equation}\label{eg1.2}
e^{- \| \beta\|_\infty t}p(t,x,y) \leq p^{\beta}(t,x,y)\le e^{\| \beta\|_\infty t}p(t,x,y) \quad\hbox{for every } t>0 \hbox{ and }  x,y\in E.
\end{equation}
Properties
\eqref{eg1.1} and \eqref{eg1.2} imply that
$\int_{E}p^{\beta}(s_{1},x,x)m(dx)<+\infty$. Thus $P^{\beta}_{t}$ is
a compact operator on $L^{2}(E,m)$ for
every $t\geq s_1$.
By Jentzch's theorem (see, for example, \cite[Theorem V.6.6]{Schaeffer}),
$-\lambda_{1}$ is a simple eigenvalue of $\mathcal{L}+\beta$ where $\mathcal{L}$ is the infinitesimal operator of $X$, and an
eigenfunction $h$ of $\mathcal{L}+\beta$ associated with
$-\lambda_{1}$ can be chosen to be positive and continuous on $E$.
Suppose
$\lambda_1 <0$.
We assume in addition that there exists $s_{2}>0$ such that
\begin{equation}
\int_{E}p(s_{2},x,x)^{2}m(dx)<+\infty.\label{eg1.4}
\end{equation}
It follows from \eqref{eg1.2} and H\"{o}lder's inequality that for
every $t>s_2$,
$P^{\beta}_{t}$ is a bounded operator from $L^{2}(E,m)$ to
$L^{4}(E,m)$.
Thus $h= e^{\lambda_1 t}P_t^\beta h\in L^{4}(E,m)$ and
so condition \eqref{llogl1} is satisfied.
Hence Theorem \ref{them2} holds.

Conditions \eqref{eg1.1} and
\eqref{eg1.4} are satisfied by a large class of Hunt processes,
which contains subordinated OU processes as special cases. By ``subordinated OU process"
we mean the process
$X_{t}=Y_{S_{t}}$,
where $Y_{t}$ is an
OU process on $\mathbb{R}^{d}$ and $S_{t}$
is a subordinator on $\mathbb{R}_{+}$ independent of $Y_{t}$. In the special case $S_{t}\equiv t$, $X_{t}$ reduces to the OU process. In general,
the sample path of $X_t$ is discontinuous.
Suppose the infinitesimal generator of $Y_{t}$ is
$$\widehat{\mathcal{L}}=\frac{1}{2}\sigma^{2}\Delta-b x\cdot\nabla\quad\mbox{ on }\mathbb{R}^{d}$$
where $\sigma,b>0$ are constants, and $S_{t}$ is a subordinator with positive drift coefficient $a>0$. As is indicated in Example 1, $Y_{t}$ is symmetric with respect to the reference measure $m(dx):=\left(\frac{b}{\pi \sigma^{2}}\right)^{d/2}\exp(-b|x|^{2}/\sigma^{2})dx$. We use $\hat{p}(t,x,y)$ to denote the transition density of $Y_{t}$ with respect to $m$. It is known that
\begin{eqnarray}
\hat{p}(t,x,y)&=&\left(1-e^{-2bt}\right)^{-d/2}\exp\left(-\frac{b}{\sigma^{2}\left(e^{2bt}-1
\right)}\left( |x|^{2}+|y|^{2}-2x\cdot y e^{bt}\right) \right) .\nonumber
\end{eqnarray}
By definition, the transition density of $X_{t}$ with respect to $m$ is given by
$$p(t,x,y)=\mathbb{E}\left[\hat{p}(S_{t},x,y)\right].$$
It is proved in \cite[Example 1.1]{RSZ} that \eqref{eg1.1} and \eqref{eg1.4} hold for a subordinated OU process.
Therefore,
Theorem \ref{them2} holds
for branching subordinated OU processes with a bounded branching rate function.

\bigskip

\noindent\textbf{Example 3. [LLN for branching diffusions on
bounded domains with branching rate given by a Kato class measure]}
Suppose $d\ge 3$, $E\subset \mathbb{R}^{d}$ is a bounded
$C^{1,1}$ domain (that is, the boundary of $E$ can be locally characterized by $C^{1,1}$ functions) and $m$ is
the Lebesgue measure on $E$. Let
$$\mathcal{L}=\frac{1}{2}\sum_{i,j=1}^{d}\partial_{i}(a_{ij}\partial_{j})$$
with $a_{ij}(x)\in C^{1}(\mathbb{R}^{d})$ for every $i,j=1,\cdots,
d$. Suppose the matrix $(a_{ij}(x))$ is symmetric and uniformly
elliptic. It is known that there exists a symmetric diffusion process $Y$ on
$\mathbb{R}^{d}$ with generator $\mathcal{L}$. Let $X$ be the killed
process of $Y$ upon $E$, i.e.,
\begin{equation} \nonumber
X_{t}=\left\{ \begin{aligned}
         Y_{t}, &\quad t<\tau_{E}, \\
  \partial,&\quad t\ge \tau_{E},
                          \end{aligned} \right.
                          \end{equation}
where $\tau_{E}:=\inf\{t>0:\ Y_{t}\not\in E\}$ and $\partial$ is a cemetery state.
Then $X$ has a transition density function
$p_{E}(t,x,y)$ which is jointly continuous in $(x,y)$ and positive
for every $t>0$. The following two-sided estimates of $p_{E}(t,x,y)$ is
established in \cite[Theorem 2.1]{Ri}, extending an earlier result of Q. Zhang.
Let
$f_{E}(t,x,y):=\left(1\wedge\frac{\delta_{E}(x)}{\sqrt{t}}\right)\left(1\wedge\frac{\delta_{E}(y)}{\sqrt{t}}\right)$,
where $\delta_{E}(x)$ denotes the distance between $x$ and the
boundary of $E$. There exist positive constants $c_{i}$,
$i=1,\cdots,4,$ such that for every $(t,x,y)\in (0,1]\times E\times
E$,
$$c_{1}f_{E}(t,x,y)t^{-d/2}e^{-c_{2}|x-y|^{2}/t}
\le p_{E}(t,x,y)\le c_{3}f_{E}(t,x,y)t^{-d/2}e^{-c_{4}|x-y|^{2}/t}.
$$
We say that a signed smooth Radon measure $\nu$ on $\mathbb{R}^{d}$
belongs to the Kato class $\mathbf{K}_{d,\alpha}$ ($\alpha\in
(0,2]$) if
\begin{equation}
\lim_{r\downarrow 0}\sup_{x\in\mathbb{R}^{d}}\int_{|x-y|\le
r}\frac{|\nu|(dy)}{|x-y|^{d-\alpha}}=0.\label{eg2.1}
\end{equation}
In fact $\mathbf{K}_{d,\alpha}$ is the Kato class
of the rotationally
symmetric
$\alpha$-stable processes on $\mathbb{R}^{d}$. We assume the branching rate measure $\mu$ is a
non-negative Radon measure in $\mathbf{K}_{d,2}$. For any
$f\in\mathcal{B}^{+}(E)$, let
$$P^{\mu}_{t}f(x):=\mathrm{E}_{x}\left[\exp( A^{\mu}_{t}) f(X_{t})\right].$$
Then $P^{\mu}_{t}$ has a transition density $p^{\mu}_{E}(t,x,y)$
which is jointly continuous in $(x,y)$ and positive for every $t>0$.
It is shown in \cite[Theorem 4.4]{Song&Kim 2} that there exist
positive constants $c_{i}$, $i=5,\cdots,8,$ such that for every
$(t,x,y)\in (0,1]\times E\times E$,
\begin{equation}
c_{5}f_{E}(t,x,y)t^{-d/2}e^{-c_{6}|x-y|^{2}/t} \le
p^{\mu}_{E}(t,x,y)\le
c_{7}f_{E}(t,x,y)t^{-d/2}e^{-c_{8}|x-y|^{2}/t}.\label{kernel}
\end{equation}
The infinitesimal generator of $P^{\mu}_{t}$ is
$(\mathcal{L}+\mu)|_{E}$ with zero Dirichlet boundary condition. It
follows from Jentzch's theorem that $-\lambda_{1}$ is a simple eigenvalue
of $(\mathcal{L}+\mu)|_{E}$ and that an eigenfunction $h$ associated
with $-\lambda_{1}$ can be chosen to be positive with
$\|h\|_{L^{2}(E,dx)}=1$. Immediately, $h$ is continuous on $E$ by the dominated convergence theorem. We assume $\lambda_{1}<0$. Recall that $E$
is bounded. Using the equation $h=e^{\lambda_{1}}P^{\mu}_{1}h$ and
the estimates in \eqref{kernel}, we get that for every $x\in E$,
\begin{equation}
c_{9}(1\wedge\delta_{E}(x))\le h(x)\le
c_{10}(1\wedge\delta_{E}(x))\label{h}
\end{equation}
for some positive constants $c_{9}$, $c_{10}$. Let
$$p^{h}_{E}(t,x,y):=\frac{e^{\lambda_{1}t}p^{\mu}_{E}(t,x,y)}{h(x)h(y)}\quad\mbox{ for }x,y\in E.$$
Immediately
condition \eqref{llogl1} holds by the boundedness of $h$ and
condition \eqref{0.2} holds by \eqref{kernel} and \eqref{h}.
Therefore both Theorem \ref{them2} and Theorem \ref{them3} hold for this example.

\bigskip
\noindent\textbf{Example 4. [LLN for branching killed
$\alpha$-stable processes with a bounded branching rate function]}
Suppose $d\ge 1$, $E=\mathbb{R}^{d}$, $m$ is the
Lebesgue measure on $\mathbb{R}^{d}$ and $\alpha\in (0,2)$. Suppose
$Y$ is a symmetric $\alpha$-stable process on $\mathbb{R}^{d}$,
and $c(x)$ is a non-negative function in $\mathbf{K}_{d,\alpha}$
( a function $q$
is said to be in $\mathbf{K}_{d,\alpha}$, if the measure
$\nu(dx):=q(x)dx$ is in $\mathbf{K}_{d,\alpha}$ where $\mathbf{K}_{d,\alpha}$ is defined in \eqref{eg2.1}). Let $X$ be the
subprocess of $Y$ such that for all
$f\in\mathcal{B}_{b}(\mathbb{R}^{d})$,
$$P_{t}f(x):=\mathrm{E}_{x}[f(X_{t})]=\mathrm{E}_{x}\left[\exp\left(
-\int_{0}^{t}c(Y_{s})ds\right) f(Y_{t})\right].$$
It is known that the infinitesimal generator of $X$ is
$\mathcal{L}=\Delta^{\alpha/2}-c(x)$, where
$\Delta^{\alpha/2}:=-(-\Delta)^{\alpha/2}$ is the generator of a
symmetric $\alpha$-stable process. Let the branching rate function
$\beta$ be a non-negative bounded function on $\mathbb{R}^{d}$. Let
$V(x):=c(x)-\beta(x)$.
Clearly,
$|V|\in
\mathbf{K}_{d,\alpha}$. For any
$f\in\mathcal{B}^{+}(\mathbb{R}^{d})$, let
$$P^{\beta}_{t}f(x):=\mathrm{E}_{x}\left[\exp\left(\int_{0}^{t}\beta(X_{s})ds\right)f(X_{t})\right].$$
Note that for every $t>0$, $P_{t}$ is bounded from
$L^{1}(\mathbb{R}^{d},dx)$ to $L^{\infty}(\mathbb{R}^{d},dx)$, and
satisfies the strong Feller property. It follows from \cite{ABM} that
for very $t>0$, $P^{\beta}_{t}$ is bounded from
$L^{p_{1}}(\mathbb{R}^{d},dx)$ to $L^{p_{2}}(\mathbb{R}^{d},dx)$ for
any $1\le p_{1}\le p_{2}\le +\infty$. Thus under our Assumption 1, the ground state
$h$ is a positive bounded continuous function on $\mathbb{R}^{d}$. The semigroup
$\{P^{\beta}_{t}:\ t\ge 0\}$ is
the Feynman-Kac semigroup with infinitesimal generator
$\mathcal{L}^{\beta}=\Delta^{\alpha/2}-V$.
Assume
in addition that $V(x)=c(x)-\beta(x)$
satisfies the following conditions:
$\liminf_{|x|\to+\infty}V^{+}(x)/\log|x|>0$ and $V^-$ has compact support.
Then by \cite[Theorem 5.5]{Kaleta}, the semigroup $P^{\beta}_{t}$ is
(AIU).
Hence both Theorem \ref{them2} and Theorem \ref{them3} are true for this example.
It is known from \cite{Kaleta} that in case of symmetric $\alpha$-stable processes,
(AIU) is weaker than (IU).
For instance, $V(x)= c(x)-\beta (x)$ with
$c(x)=\log|x|1_{\{|x|\ge R\}}$ and $\beta$ has compact support in $ B(0,R)$
for some $R\geq 1$ is such an example.

\medskip

The rest of the paper is organized as follows. In Section 2, we review
some facts on  Girsanov transform and $h$-transforms
in the context of symmetric Markov processes, and prove Proposition \ref{P:1.8}.
Spine construction of branching processes is recalled in Section 3.
We then present proof for the $L\log L$ criteria, Theorem \ref{them1}, in Section 4.
Weak law of large numbers, Theorem \ref{them2}, is proved in Section 5, while Theorem
\ref{them3} on the strong law of large numbers will be proved in Section 6.
The lower case constants
$c_{1},c_{2},\cdots,$ will denote the generic constants used in this
article, whose exact values are not important, and can change from
one appearance to another.

\section{Preliminary}

Recall $h\in \FF^\mu$ is the minimizer in Assumption 1.
Since $h\in \mathcal{F}$, by Fukushima's decomposition, we have for
q.e. $x\in E$, $\Pi_{x}$-a.s.
$$
h(X_{t})-h(X_{0})=M^{h}_{t}+N^{h}_{t},\quad    t\ge 0,
$$
where $M^{h}$ is a martingale additive functional of $X$ having
finite energy and $N^{h}_{t}$ is a continuous additive functional of
$X$ having zero energy. It follows from \eqref{1.3} and \cite[Theorem 5.4.2]{FOT} that $N^{h}_{t}$ is of
bounded variation, and
$$N^{h}_{t}=-\lambda_{1}\int_{0}^{t}h(X_{s})ds-\int_{0}^{t}h(X_{s})dA^{(Q-1)\mu}_{s},\quad\forall t\ge 0.$$
Following \cite[Section 2]{CFTYZ} (see also \cite[Section 2]{Chen&Shiozawa}),  we define
a local martingale on the random time interval
$[0,\zeta_{p})$ by
\begin{equation} \label{(1)}
M_{t}:=\int_{0}^{t}\frac{1}{h(X_{s-})}dM^{h}_{s},
\quad  t\in [0,\zeta_{p}),
\end{equation}
where $\zeta_{p}$ is the predictable part of the life time $\zeta$
of $X$. Then the solution $R_{t}$ of the stochastic differential
equation
\begin{equation}
R_{t}=1+\int_{0}^{t}R_{s-}dM_{s},  \quad
 t\in [0,\zeta_{p}),\label{(2)}
\end{equation}
is a positive local martingale on $[0,\zeta_{p})$ and hence a
supermartingale. As a result, the formula
$$d\Pi^{h}_{x}=R_{t}d\Pi_{x}\quad \mbox{on }\mathcal{H}_{t}\cap\{t<\zeta\}\quad\mbox{ for }x\in E$$
uniquely determines a family of subprobability measures
$\{\Pi^{h}_{x}:x\in E\}$ on $(\Omega,\mathcal{H})$. We denote $X$
under $\{\Pi^{h}_{x}:x\in E\}$ by $X^{h}$, that is
$$\Pi^{h}_{x}\left[f(X^{h}_{t})\right]=\Pi_{x}\left[R_{t}f(X_{t}):t<\zeta\right]$$
for any $t\ge 0$ and $f\in\mathcal{B}^{+}(E)$. It follows from
\cite[Theorem 2.6]{CFTYZ} that the process $X^{h}$ is an irreducible recurrent
$\widetilde m$-symmetric right Markov process, where $\widetilde{m}(dy):=h(y)^{2}m(dy)$.
 Note that by \eqref{(1)},
\eqref{(2)} and Dol\'{e}an-Dade's formula,
\begin{equation}\label{(3)}
R_{t}=\exp\big( M_{t}-\frac{1}{2}\langle M^{c}\rangle_{t}\big) \prod_{0<s\le
t}\frac{h(X_{s})}{h(X_{s-})}\exp\left( 1-\frac{h(X_{s})}{h(X_{s-})}\right),
\quad t\in [0, \zeta_p),
\end{equation}
where $M^{c}$ is the continuous martingale part of $M$. Applying
Ito's formula to $\log h(X_{t})$, we obtain that for q.e. $x\in E$,
$\Pi_{x}$-a.s. on  $[0,\zeta)$,
\begin{equation}
\log h(X_{t})-\log h(X_{0})=M_{t}-\frac{1}{2}\langle
M^{c}\rangle_{t}+\sum_{s\le t}
\left(\log\frac{h(X_{s})}{h(X_{s-})}-\frac{h(X_{s})-h(X_{s-})}{h(X_{s-})}\right)-\lambda_{1}t-A^{(Q-1)\mu}_{t}.\label{(4)}
\end{equation}
By \eqref{(3)} and \eqref{(4)}, we get
$$
R_{t}=\exp\big( \lambda_{1}t+A^{(Q-1)\mu}_{t}\big)
\frac{h(X_{t})}{h(X_{0})}\quad\mbox{on }[0,\zeta).
$$
Therefore for any $f\in\mathcal{B}^{+}(E)$,
\begin{eqnarray} \label{1.7}
\Pi^{h}_{x}\left[f(X^{h}_{t})\right]
=\frac{e^{\lambda_{1}t}}{h(x)}\Pi_{x}\left[e^{A^{(Q-1)\mu}_{t}}h(X_{t})f(X_{t})\right]
 = \frac{e^{\lambda_{1}t}}{h(x)}P^{(Q-1)\mu}_{t}(hf)(x)
 =P^{h}_{t}f(x) .
\end{eqnarray}
Let $(\mathcal{E}^{h},\mathcal{F}^{h})$ be the symmetric Dirichlet
form on $L^{2}(E,\widetilde{m})$ generated by $X^{h}$.
Then \eqref{1.7} says that the transition semigroup of $X^h$ is exactly the semigroup
$\{P^{h}_{t}: t\geq 0\} $
obtained from $P^{(Q-1)\mu}_{t}$ through Doob's $h$-transform.
Consequently,  $f\in \mathcal{F}^{h}$ if and only if $fh\in\mathcal{F}^{\mu}$, and
$$
\mathcal{E}^{h}(f,f)=\mathcal{E}^{(Q-1)\mu}(fh,fh)-\lambda_{1}\int_{E}f(x)^{2}h(x)^{2}m(dx).
$$
In other words, $\Phi^{h}:f\mapsto fh$ is an isometry from $(\mathcal{E}^{h},\mathcal{F}^{h})$ onto $(\mathcal{E}^{(Q-1)\mu+\lambda_{1}m},\mathcal{F}^{\mu})$ and from $L^{2}(E,\widetilde{m})$ onto $L^{2}(E,m)$.
Let
$\sigma(\mathcal{E}^{h})$ denote the spectrum of the positive
definite self-adjoint operator associated with $\mathcal{E}^{h}$.
We know from  \cite[Theorem 2.6]{CFTYZ} that the constant function $1$ belongs to
$\mathcal{F}^{h}$, and $\mathcal{E}^{h}(1,1)=0$.
Hence $0\in \sigma(\mathcal{E}^{h})$ is a simple
eigenvalue and $1$ is the corresponding eigenfunction.
Therefore
$$
\lambda^{h}_{1}:=\inf\left\{\mathcal{E}^{h}(u,u):\ u\in\mathcal{F}^{h} \hbox{ with }  \int_{E}u(x)^{2}\widetilde{m}(dx)=1\right\}=0.
$$
Let $\lambda^{h}_{2}$ be the second bottom of $\sigma(\mathcal{E}^{h})$, i.e.
$$
\lambda^{h}_{2}:=\inf\left\{\mathcal{E}^{h}(u,u):\ u\in\mathcal{F}^{h} \hbox{ with }  \int_{E}u(x)\widetilde{m}(dx)=0
\hbox{ and }  \int_{E}u(x)^{2}\widetilde{m}(dx)=1\right\}.
$$
In view of the isometry $\Phi^h$,  we have $\lambda^{h}_{2}= \lambda_2 -\lambda_1$.
So Assumption 2 is equivalent to assuming $\lambda^h_2>0$.

Define
\begin{equation}\label{def-a}
\widetilde{a}_{t}(x):=p^{h}(t,x,x)
\quad \mbox{ for }  t>0 \mbox{ and } x\in E.\end{equation}
Note that by \eqref{p1} for any $x\in E$ and $t,s\ge 0$,
\begin{eqnarray}
|P^{h}_{t+s}\varphi(x)|&=&|P^{h}_{s}P^{h}_{t}\varphi(x)|\nonumber\\
&\le&\int_{E}p^{h}(s,x,y)|P^{h}_{t}\varphi(y)|\widetilde{m}(dy)\nonumber\\
&\le&\left(\int_{E}p^{h}(s,x,y)^{2}\widetilde{m}(dy)\right)^{1/2}\|P^{h}_{t}\varphi\|_{L^{2}(E,\widetilde{m})}\nonumber\\
&\le&\widetilde{a}_{2s}(x)^{1/2}e^{-\lambda_h t}\|\varphi\|_{L^{2}(E,\widetilde{m})}.\label{p2}
\end{eqnarray}
We use $\langle f,g\rangle_{L^{2}(\widetilde{m})}$ to denote $\int_{E}f(x)g(x)\widetilde{m}(dx)$.
For every $g\in L^{2}(E,\widetilde{m})$ we can decompose  $g$ as
$g(x)=\langle g,1\rangle_{L^2(\widetilde m)}  +\varphi(x)$
with $\varphi(x)$ satisfying $\int_{E}\varphi(x)\widetilde{m}(dx)=0$. Then by \eqref{p2}, for any $t>s\ge 0$ and any $x\in E$,
\begin{eqnarray}
  |P^{h}_{t}g(x)-\langle g,1\rangle_{L^2(\widetilde m)}  |
&=&|P^{h}_{t}\varphi(x)|
\leq e^{-\lambda_h (t-s)}\widetilde{a}_{2s}(x)^{1/2}\|\varphi\|_{L^{2}(E,\widetilde{m})}\nonumber\\
&\le & e^{-\lambda_h (t-s)}\widetilde{a}_{2s}(x)^{1/2}\|g\|_{L^{2}(E,\widetilde{m})}.\label{p3}
\end{eqnarray}

\medskip

\noindent{\textit{Proof of Proposition \ref{P:1.8}.}
We only need to  prove that
\eqref{0.2} implies  \eqref{1.17}.
For any $f\in L^{2}(E,\widetilde{m})$, define
$c_{f}:=\int_{E}f(x)\widetilde{m}(dx)$. Immediately, we have for any $t>0$,
$\int_{E}(f-c_{f})\widetilde{m}(dx)=0$, and $f-c_{f}\in
L^{2}(E,\widetilde{m})$. By \eqref{p1}, we have
\begin{equation}
\|P^{h}_{t}f-c_{f}\|_{L^{2}(E,\widetilde{m})}\le
e^{-\lambda_h t}\|f\|_{L^{2}(E,\widetilde{m})}.\label{poincare}
\end{equation}
Let $t_{1}>0$ be the constant in \eqref{0.2}.
By the semigroup property, for any $t\ge t_{1}/2$ and $x\in E$,
\begin{eqnarray}
p^{h}(t,x,y)&=&\int_{E}p^{h}(t_{1}/2,x,z)p^{h}(t-t_{1}/2,y,z)\widetilde{m}(dz)\nonumber\\
&=&P^{h}_{t-t_{1}/2}f_{x}(y),\label{1.14}
\end{eqnarray}
where $f_{x}(\cdot):=p^{h}(t_{1}/2,x,\cdot)$.  Note that by
H\"{o}lder's inequality,
\begin{eqnarray}\label{0.4}
|p^{h}(t,x,y)-1|&=&|\int_{E}p^{h}(t/2,x,z)p^{h}(t/2,z,y)\widetilde{m}(dz)-1|\nonumber\\
&=&|\int_{E}(p^{h}(t/2,x,z)-1)(p^{h}(t/2,z,y)-1)\widetilde{m}(dz)|\nonumber\\
&\le&\left(\int_{E}(p^{h}(t/2,x,z)-1)^{2}\widetilde{m}(dz)\right)^{1/2}
\left(\int_{E}(p^{h}(t/2,z,y)-1)^{2}\widetilde{m}(dz)\right)^{1/2}
\end{eqnarray}
Note that
$c_{f_{x}}=\int_{E}p^{h}(t_{1}/2,x,y)\widetilde{m}(dy)=1$ and
$\int_{E}f_{x}^{2}(y)\widetilde{m}(dy)=\widetilde{a}_{t_{1}}(x)$. By \eqref{1.14} and \eqref{poincare},
for any $t>t_{1},$
\begin{eqnarray}
\left(\int_{E}(p^{h}(t/2,x,z)-1)^{2}\widetilde{m}(dz)\right)^{1/2}&=&\|P^{h}_{(t-t_{1})/2}f_{x}-c_{f_{x}}\|_{L^{2}(E,\widetilde{m})}\nonumber\\
&\le&e^{-\lambda_h (t-t_{1})/2}\|f_{x}\|_{L^{2}(E,\widetilde{m})}\nonumber\\
&=&e^{-\lambda_h (t-t_{1})/2}\widetilde{a}_{t_{1}}(x)^{1/2}.\label{0.5}
\end{eqnarray}
Similarly, for any $t>t_{1},$
\begin{equation}
\left(\int_{E}(p^{h}(t/2,z,y)-1)^{2}\widetilde{m}(dz)\right)^{1/2}\le
e^{-\lambda_h (t-t_{1})/2}\widetilde{a}_{t_{1}}(y)^{1/2}. \label{0.6}
\end{equation}
Combining \eqref{0.4}, \eqref{0.5} and \eqref{0.6},
we have for any $t>t_{1}$,
\begin{equation}\label{new-eq2}|p^{h}(t,x,y)-1|\le e^{-\lambda_h
(t-t_{1})}\widetilde{a}_{t_{1}}(x)^{1/2}\widetilde{a}_{t_{1}}(y)^{1/2}.\end{equation}
 Therefore \eqref{0.2} implies
that for any $t>t_{1}$,  \eqref{1.17} holds.
 \qed

 \section{Spine construction}

 To establish the $L^{1}$ convergence of the martingale $M_{t}$,
 we apply the ``spine" and change of measure techniques presented
  in \cite{HH} for branching diffusions
  to our branching Hunt processes. The notation used here is closely related to that used in \cite{HH}. It is known that the family structure of the particles in a branching Markov process can be well expressed by marked Galton-Watson trees
  (see, for example, \cite{HH, LRSb} and the references therein).
   Let $\mathcal{T}$ denote the space of all marked G-W trees. For a fixed $\tau\in \mathcal{T}$, all particles in $\tau$ are labeled according to the Ulam-Harris convention, for example, $\emptyset231$ or $231$ is the first child of the third child of the second child of the initial ancestor $\emptyset$. Besides, each particle $u\in\tau$ has a mark $(X_{u},\sigma_{u},A_{u})$ where $X_{u}:[b_{u},\zeta_{u})\to E_{\partial}$ is the spatial location of $u$ during its life
  span ($b_{u}$ is its birth time and $\zeta_{u}$ its fission time),
 $\sigma_{u}=\zeta_{u}-b_{u}$
 is the length of its life span,
 and $A_{u}$ denotes the number of its offspring. We use $u\prec v$ to mean that $u$ is an ancestor of $v$.

Since in this paper every particle is assumed to give birth to at least one child,
for each tree $\tau$, we can choose a distinguished genealogical
 path of descent from the initial ancestor $\emptyset$. Such a line is called a \textit{spine} and denoted as $\xi=\{\xi_{0}=\emptyset,\xi_{1},\xi_{2},\cdots\}$, where $\xi_{i}$ is the label of the $i$th
 spine node.
Define $\mbox{node}_{t}(\xi):=u$ if $u\in\xi$ and is alive at time $t$.
 We shall use $\{\widetilde X_{t}:t\ge 0\}$ and $\{n_{t}:t\ge 0\}$ respectively to denote the spatial path and the counting process of fission times
 along the spine. Let $\widetilde{\mathcal{T}}$ denote the space of G-W trees with a distinguished spine. We introduce some fundamental filtrations that encapsulate different knowledge:
 $$
 \mathcal{F}_{t}:=\sigma\left\{ (u,X_{u},\sigma_{u},A_{u}),  \zeta_{u}\le t;\ (u,X_{u}(s), s\in[b_{u},t]), t\in [b_{u},\zeta_{u})\right\},\quad \mathcal{F}_{\infty}:=\bigvee_{t\ge 0}\mathcal{F}_{t}
 := \sigma \{ \mathcal{F}_t; t\geq 0\};
 $$
 $$
 \widetilde{\mathcal{F}}_{t}:=\sigma\left\{\mathcal{F}_{t},\mbox{node}_{t}(\xi)\right\},
 \quad \widetilde{\mathcal{F}}_{\infty}=\bigvee_{t\ge 0}\widetilde{\mathcal{F}}_{t};$$
 $$\mathcal{G}_{t}:=\sigma \{\widetilde{X}_{s}:s\le t\},\quad \mathcal{G}_{\infty}:=\bigvee_{t\ge 0}\mathcal{G}_{t};$$
 $$\widetilde{\mathcal{G}}_{t}:=\sigma\left\{\mathcal{G}_{t},(\mbox{node}_{s}(\xi):s\le t),(A_{u}:\ u\prec \mbox{node}_{t}(\xi))\right\},\quad \widetilde{\mathcal{G}}_{\infty}:=\bigvee_{t\ge 0}\widetilde{\mathcal{G}}_{t}.
 $$
We assume $\mathbb{P}_{x}$ is the measure on $( \widetilde{\mathcal{T}},
\mathcal{F}_{\infty})$
such that the filtered probability space $\left( \widetilde{\mathcal{T}},
\mathcal{F}_{\infty},\left(\mathcal{F}_{t}\right)_{t\ge 0},\mathbb{P}_{x}\right)$ is the canonical model for the branching Hunt process $\mathbb{X}$ described in the introduction. We know from \cite{HH} that the measure $\mathbb{P}_{x}$ on $( \widetilde{\mathcal{T}},
\mathcal{F}_{\infty})$ can be extended to the probability measure $\widetilde{\mathbb{P}}_{x}$ on $( \widetilde{\mathcal{T}},
\widetilde{\mathcal{F}}_{\infty})$ such that the $n$th spine node is uniformly chosen from the children of the $(n-1)$th spine node.

For every $t\ge 0$, as in \cite{HH}, we define
$$
\eta (t)
:=\exp\left(\lambda_{1}t+A^{(Q-1)\mu}_{t} \right)
\frac{h(\widetilde{X}_{t})}{h(x)},\quad Z(t):=e^{\lambda_{1}t}\frac{\mathbb{X}_{t}(h)}{h(x)}, \quad\mbox{ and  }\quad\widetilde{\eta}(t):=e^{\lambda_{1}t}\frac{h(\widetilde X_{t})}{h(x)}\prod_{v\prec
\mbox{node}_{t}(\xi)}A_{v}.$$
Then $\eta(t)$, $Z(t)$ and $\widetilde{\eta}(t)$ are positive $\widetilde{\mathbb{P}}_{x}$-martingales with respect to the filtrations $\mathcal{G}_{t}$, $\mathcal{F}_{t}$ and $\widetilde{\mathcal{F}}_{t}$, respectively.
Moreover, both $\eta(t)$ and $Z(t)$ are projections of $\widetilde{\eta}(t)$ onto their filtrations:
$$Z(t)=\widetilde{\mathbb{P}}_{x}\left[\widetilde{\eta}(t)|\mathcal{F}_{t}\right],\quad \eta(t)=\widetilde{\mathbb{P}}_{x}\left[\widetilde{\eta}(t)|\mathcal{G}_{t}\right]\quad\mbox{ for }t\ge 0.$$
We call  $\eta(t)$ the \textit{single-particle martingale} and $Z(t)$  the \textit{branching-particle martingale}.
As
in \cite{HH}, we define  a new probability measure $\widetilde{\mathbb{Q}}_{x}$ by setting
 \begin{equation}
 d\widetilde{\mathbb{Q}}_{x}|_{\widetilde{\mathcal{F}}_{t}}=\widetilde{\eta}(t)
 d\widetilde{\mathbb{P}}_{x}|_{\widetilde{\mathcal{F}}_{t}}\quad\mbox{ for }t\ge 0.\label{2.1}
 \end{equation}
This implies that
\begin{equation}
 d\widetilde{\mathbb{Q}}_{x}|_{\mathcal{F}_{t}}=Z(t)
 d\widetilde{\mathbb{P}}_{x}|_{\mathcal{F}_{t}}\quad\mbox{ for }t\ge 0.
 \end{equation}
The influence of the measure change \eqref{2.1} lies in the following three aspects: Firstly, under $\widetilde{\mathbb{Q}}_{x}$ the motion of the spine is biased by the
martingale $\eta(t)$. Secondly, the branching events along
the spine occur at an accelerated rate. Finally, the number of children of the spine
nodes is size-biased distributed, that is, for every spine node $v$ located at $x$, $A_{v}$ is distributed
according to the probability mass function
 $\{\widehat{p}_{k}(x):=kp_{k}(x)/Q(x):k=0,1,\cdots\}$, while other
 (non-spine) nodes, once born, remain unaffected. More specifically, under measure
$\widetilde{\mathbb{Q}}_{x}$,
\begin{enumerate}
\item[(i)] The spine's spatial process $\widetilde{X}$ moves in $E$ as a copy of $(X^{h},\Pi^{h}_{x})$;

\item[(ii)] The branching events along
the spine occur at an accelerated rate $\widetilde{\mu}(dx):=Q(x)h(x)^{2}\mu(dx)$.
This implies that given $\mathcal{G}_{\infty}$, for every $t>0$, the number of fission times $n_{t}$ is a Poisson random variable with parameter
$ \widetilde{A}_{t}$,
where $\widetilde{A}_{t}$ is a PCAF of $\widetilde{X}$
having  Revuz measure $\widetilde{\mu}$.
To emphasize this dependence,
we also write $A^{\widetilde{\mu}}_{t}$ for $\widetilde{A}_{t}$.

\item[(iii)] At the fission time of node $v$ in the spine, the
single spine particle is replaced by $A_{v}$ children, with $A_{v}$
being distributed according to the size-biased distribution
$\{\widehat{p}_{k}(\widetilde{X}_{\zeta_{v}-}):k=0,1,\cdots\}$. Each child is
selected
to be
the next spine node with equal probabilities.

\item[(iv)] Each of the remaining $A_{v}-1$ children gives rise to independent subtrees, which are not part of the spine and evolve as independent processes determined by the measure $\mathbb{P}$ shifted to their point and time of creation.
\end{enumerate}
For more details on
martingale changes of measures for branching Hunt processes, see \cite{HH} and \cite{LRSb}.

\section{$L\log L$ criteria}

In this section, we prove Theorem \ref{them1}.
It follows from \cite[Theorem 4.1.1]{CF} that if $\nu$
is a smooth measure of $X^{h}$, then for every $g\in
\mathcal{B}^{+}(E)$ and $t\ge 0$,
\begin{equation}
\Pi^{h}_{\widetilde{m}}\left[\int_{0}^{t}g(X^{h}_{s})d
A^{\nu}_{s}\right]=\int_{0}^{t}\int_{E}g(y)\nu(dy)ds.\label{1.8}
\end{equation}
Since $t$ is arbitrary, the monotone class theorem implies that for
any $f(s,x)=l(s)g(x)$ with $l\in\mathcal{B}^{+}[0,+\infty)$ and
$g\in\mathcal{B}^{+}(E)$,
\begin{equation}
\Pi^{h}_{\widetilde{m}}\left[\int_{0}^{t}f(s,X^{h}_{s})d
A^{\nu}_{s}\right]=\int_{0}^{t}\int_{E}f(s,y)\nu(dy)ds.\label{equa1}
\end{equation}
Note that for every $f\in \mathcal{B}^{+}([0,+\infty)\times E)$,
there exists a sequence of functions $f_{n}(s,x)=l_{n}(s)g_{n}(x)$
with $l_{n}\in \mathcal{B}^{+}[0,+\infty)$ and $g_{n}\in
\mathcal{B}^{+}(E)$ such that $f_{n}(s,x)$ converges increasingly to
$f(s,x)$ for all $s$ and $x$. Thus by the monotone convergence
theorem, \eqref{equa1} holds for all $f\in
\mathcal{B}^{+}([0,+\infty)\times E)$.

\bigskip

\noindent\textit{Proof of Theorem \ref{them1}.} Since $M_{t}$ is a
positive martingale, it suffices to prove that
$\mathbb{E}_{x}M_{\infty}=h(x)$ for all $x\in E$.
By \cite[Theorem 5.3.3]{Durrett} (or \cite[Lemma 3.1]{LRSb}),
it suffices to show
that if \eqref{llogl} holds, then
\begin{equation} \label{3.1}
 \widetilde{\mathbb{Q}}_{x} \left(\limsup_{\mathbb{N}\ni n\to+\infty} M_{n\sigma}<+\infty \right)=1   \quad\hbox{for every } \sigma>0 \hbox{ and }  x\in E.
 \end{equation}
Recall that $\widetilde{\mathcal{G}}_{\infty}$
contains all the information about the spine.
We have the following \textit{spine decomposition} for $M_{t}$,
 \begin{eqnarray}
  \widetilde{\mathbb{Q}}_{x}\left[M_{t} \,\big|\, \widetilde{\mathcal{G}}_{\infty}\right
  ]&=&\widetilde{\mathbb{Q}}_{x}\left[e^{\lambda_{1}t}\sum_{u\in{\cal Z}_{t}}h(X_{u}(t)) \,\Big|\, \widetilde{\mathcal{G}}_{\infty}\right]\nonumber\\
  &=&e^{\lambda_{1}t}h(\widetilde X_{t})+\sum_{u\prec \xi_{n_{t}}}(A_{u}-1)e^{\lambda_{1}\zeta_{u}}E_{\widetilde X_{\zeta_{u}}}\left[e^{\lambda_{1}(t-\zeta_{u})}\mathbb{X}_{t-\zeta_{u}}(h)\right]\nonumber\\
  &=&e^{\lambda_{1}t}h(\widetilde X_{t})+\sum_{u\prec \xi_{n_{t}}}(A_{u}-1)e^{\lambda_{1}\zeta_{u}}h(\widetilde X_{\zeta_{u}}).\label{spinedecomp}
  \end{eqnarray}
Let $\mathcal{G}_{t}:=\sigma\{\widetilde X_{s}:\ s\le t\}$ and
$\mathcal{G}_{0}$ be the trivial $\sigma$-field. It follows from the
second Borel-Cantelli lemma (see, for example,
\cite[Section 5.4]{Durrett}) and the Markov property that
for any $\sigma >0$ and $M\geq 1$,
\begin{eqnarray}
\widetilde{\mathbb{Q}}_{\widetilde{m}}\left(e^{\lambda_{1}n\sigma}h(\widetilde X_{n\sigma})>M\mbox{ i.o.}\right)
&=&\widetilde{\mathbb{Q}}_{\widetilde{m}}\left(\sum_{n=1}^{+\infty}\widetilde{\mathbb{Q}}_{\widetilde{m}}
\left(e^{\lambda_{1}n\sigma}h(\widetilde X_{n\sigma})>M~ \Big|~\mathcal{G}_{(n-1)\sigma}\right)=+\infty\right)\nonumber\\
&=&\Pi^{h}_{\widetilde{m}}\left(\sum_{n=1}^{+\infty}\Pi^{h}_{\widetilde X_{(n-1)\sigma}}\left(e^{\lambda_{1}n\sigma}h(\widetilde X_{\sigma})>M\right)=+\infty\right).\label{3.22}
\end{eqnarray}
Recall that $\widetilde{m}$ is the invariant probability measure of
$(\widetilde X,\Pi^{h}_{x})$. Thus by Fubini's theorem and Markov property,
we have
\begin{eqnarray}
\Pi^{h}_{\widetilde{m}}\left[\sum_{n=1}^{+\infty}\Pi^{h}_{\widetilde X_{(n-1)\sigma}}\left(e^{\lambda_{1}n\sigma}h(\widetilde X_{\sigma})>M\right)\right]\label{3.2}
&=&\sum_{n=1}^{+\infty}\Pi^{h}_{\widetilde{m}}\left(e^{\lambda_{1}n\sigma}h(\widetilde X_{n\sigma})>M\right)\nonumber\\
&=&\sum_{n=1}^{+\infty}\int_{E}1_{\{e^{\lambda_{1}n\sigma}h(y)>M\}}\widetilde{m}(dy)\nonumber\\
&=&\sum_{n=1}^{+\infty}
\widetilde m \Big( \frac{\log^{+}h(y)-\log M}{-\lambda_{1}}>n\sigma \Big)\nonumber\\
&\leq & (-\lambda_1)^{-1} \int_E \log^+ h(y) \widetilde m (dy)<\infty. \nonumber
\end{eqnarray}
Therefore by \eqref{3.22} we have
$\widetilde{\mathbb{Q}}_{\widetilde{m}}\left(e^{\lambda_{1}n\sigma}h(\widetilde X_{n\sigma})>M\mbox{
i.o.}\right)=0.$ Consequently,
$$\widetilde{\mathbb{Q}}_{m}\left(\limsup_{\mathbb{N}\ni
n\to+\infty}e^{\lambda_{1}n\sigma}h(\widetilde X_{n\sigma})<+\infty\right)=1.$$
It is easy to check that the function
$x\mapsto\widetilde{\mathbb{Q}}_{x}\left(\limsup_{\mathbb{N}\ni
n\to+\infty}e^{\lambda_{1}n\sigma}h(\widetilde X_{n\sigma})<+\infty\right)$
is an invariant function for $(\widetilde X,\Pi^{h}_{\widetilde{m}})$.
Recall that $\widetilde{X}$ has a transition density function with respect to $\widetilde{m}$.
By \cite[Theorem A.2.17]{CF},
$$
\widetilde{\mathbb{Q}}_{x}\left(\limsup_{\mathbb{N}\ni
n\to+\infty}e^{\lambda_{1}n\sigma}h(\widetilde X_{n\sigma})<+\infty\right)=1
 \quad \mbox{for every } x\in E.
$$
Suppose $\eps \in (0,-\lambda_{1})$. For simplicity we
use $\zeta_{i}$ and $A_{i}$ to denote respectively the
fission time
and offspring number of the $i$th spine node.
\begin{eqnarray}
 \sum_{i=1}^{+\infty}e^{\lambda_{1}\zeta_{i}}A_{i}h(\widetilde X_{\zeta_{i}})&=& \sum_{i=1}^{+\infty}e^{\lambda_{1}\zeta_{i}}A_{i}h(\widetilde X_{\zeta_{i}})1_{\{A_{i}h(\widetilde X_{\zeta_{i}})\le e^{\eps  \zeta_{i}}\}}
 +\sum_{i=1}^{+\infty}e^{\lambda_{1}\zeta_{i}}A_{i}h(\widetilde X_{\zeta_{i}})1_{\{A_{i}h(\widetilde X_{\zeta_{i}})>e^{\eps  \zeta_{i}}\}}\nonumber\\
 &=:&I+II.
\end{eqnarray}
Recall that $\mathcal{G}_{\infty}$
contains all the
information of spine's motion
$\{\widetilde{X}_{t}: t\ge 0\}$
in $E$.
For any set $B\in
\mathcal{B}[0,+\infty)\times \mathcal{B}(\mathbb{Z}_{+})$, define
$N(B):=\sharp\{i\ge 1:\ (\zeta_{i},A_{i})\in B\}$. Then conditioned
on $\mathcal{G}_{\infty}$, $N$ is a Poisson random measure on
$[0,+\infty)\times \mathbb{Z}_{+}$ with intensity
$dA^{\widetilde{\mu}}_{s}\sum_{k\in\mathbb{Z}_{+}}\hat{p}_{k}(\widetilde X_{s})\delta_{k}(dy)$, where $\widetilde{\mu}(dx)=Q(x)h(x)^{2}\mu(dx)$.
We have
 $$\widetilde{\mathbb{Q}}_{x}\left[\sum_{i=1}^{+\infty}1_{\{A_{i}h(\widetilde X_{\zeta_{i}})>e^{\eps  \zeta_{i}}\}}\right]=\Pi^{h}_{x}\left[\int_{0}^{+\infty}
 \sum_{k=0}^{+\infty}\hat{p}_{k}(\xi_{s})1_{\{kh(\widetilde X_{s})>e^{\eps  s}\}}d A^{\widetilde{\mu}}_{s}\right].$$
 Note that by \eqref{equa1} and our assumption \eqref{llogl},
 \begin{eqnarray}
 \widetilde{\mathbb{Q}}_{\widetilde{m}}\left[\sum_{i=1}^{+\infty}1_{\{A_{i}h(\widetilde X_{\zeta_{i}})>e^{\eps  \zeta_{i}}\}}\right]&=&\Pi^{h}_{\widetilde{m}}\left[\int_{0}^{+\infty}
 \sum_{k=0}^{+\infty}\hat{p}_{k}(\xi_{s})1_{\{kh(\widetilde X_{s})>e^{\eps  s}\}}d A^{\widetilde{\mu}}_{s}\right]\nonumber\\
 &=&\int_{0}^{+\infty}\int_{E}\sum_{k=0}^{+\infty}\hat{p}_{k}(y)1_{\{\log^{+}(k h(y))>\eps  s\}}\widetilde{\mu}(dy)ds\nonumber\\
 &=&\eps ^{-1}\int_{E}\sum_{k=0}^{+\infty}kp_{k}(y)\log^{+}(k h(y))h(y)^{2}\mu(dy)<+\infty.\nonumber
 \end{eqnarray}
Thus
$$\widetilde{\mathbb{Q}}_{x}\left(\sum_{i=1}^{+\infty}1_{\{A_{i}h(\widetilde X_{\zeta_{i}})>e^{\eps  \zeta_{i}}\}}<+\infty\right)=1\quad \hbox{for }
\widetilde{m}
\hbox{-a.e. }
x\in E.
$$
This implies that $II$ is
the sum of a finite many terms and so
$\widetilde{\mathbb{Q}}_{x}(II<+\infty)=1$ for $\widetilde{m}$-a.e.
$x\in E$. On the other hand, since $Q(x)$ is bounded on $E$, we have
\begin{eqnarray}
 \widetilde{\mathbb{Q}}_{\widetilde{m}}[I]&=&\Pi^{h}_{\widetilde{m}}\left[\int_{0}^{+\infty}e^{\lambda_{1}s}
 \sum_{k=0}^{+\infty}\hat{p}_{k}(\widetilde X_{s})k h(\widetilde X_{s})1_{\{k h(\xi_{s})\le e^{\eps  s}\}}d A^{\widetilde{\mu}}_{s}\right]\nonumber\\
 &=&\int_{0}^{+\infty}\int_{E}e^{\lambda_{1}s}\sum_{k=0}^{+\infty}\hat{p}_{k}(y)k h(y)1_{\{k h(y)\le e^{\eps  s}\}}\widetilde{\mu}(dy)ds\nonumber\\
 &\le&\int_{0}^{+\infty}\int_{E} e^{(\lambda_{1}+\eps )s}\sum_{k=0}^{+\infty}\hat{p}_{k}(y)Q(y)h(y)^{2}\mu(dy)ds\nonumber\\
 &\le&\|Q\|_{\infty}\int_{E}h(y)^{2}\mu(dy)<+\infty.\nonumber
 \end{eqnarray}
Thus we have $\widetilde{\mathbb{Q}}_{x}(I<+\infty)=1$ for
$\widetilde{m}$-a.e. $x\in E$. Now we have proved that
$$\widetilde{\mathbb{Q}}_{x}\left(\sum_{i=1}^{+\infty}e^{\lambda_{1}\zeta_{i}}A_{i}h(\widetilde X_{\zeta_{i}})<+\infty\right)=1\quad \hbox{for }  \widetilde{m}\mbox{-a.e. }x\in E.$$
It is easy to check that the function
$x\mapsto\widetilde{\mathbb{Q}}_{x}\left(\sum_{i=1}^{+\infty}e^{\lambda_{1}\zeta_{i}}A_{i}h(\widetilde X_{\zeta_{i}})<+\infty\right)$
is an invariant function. Thus by
\cite[Theorem A.2.17]{CF},
we get
$\widetilde{\mathbb{Q}}_{x}\left(\sum_{i=1}^{+\infty}e^{\lambda_{1}\zeta_{i}}A_{i}h(\widetilde X_{\zeta_{i}})<+\infty\right)
=1$ for every $x\in E$.
By \eqref{spinedecomp} we have
$$\widetilde{\mathbb{Q}}_{x}\left(\limsup_{\mathbb{N}\ni n\to +\infty}\widetilde{\mathbb{Q}}_{x}\left[M_{n\sigma}\,\big|\,\widetilde{\mathcal{G}}_{\infty}\right]
<+\infty\right)=1 \quad \hbox{for every } x\in E.
$$
By Fatou's lemma, we get
$\widetilde{\mathbb{Q}}_{x}\left(\liminf_{\mathbb{N}\ni n\to+\infty}M_{n\sigma}<+\infty\right)=1$.
Note that $M_{n\sigma}^{-1}$ is a positive
$\widetilde{\mathbb{Q}}_{x}$-martingale with respect to
$\mathcal{F}_{n\sigma}$,
so it converges almost surely as $n\to \infty$. It follows then
$$\widetilde{\mathbb{Q}}_{x}\left(\limsup_{\mathbb{N}\ni n\to+\infty}M_{n\sigma}=\liminf_{\mathbb{N}\ni n\to+\infty}M_{n\sigma}<+\infty\right)=1.$$\qed

\section{Weak law of large numbers}

In this section, we present a proof for Theorem \ref{them2}}.

\begin{lemma}\label{lem2.1}
If Assumption 1 and \eqref{llogl} hold, then for every $t\ge 0$
and $\phi\in \mathcal{B}^{+}_{b}(E)$,
\begin{equation}\label{eq2.1}
\int_{E}\mathbb{E}_{x}\left[\mathbb{X}_{t}(\phi
h)\log^{+}\mathbb{X}_{t}(\phi h)\right]h(x)m(dx)<+\infty.
\end{equation}
\end{lemma}

\proof First we note that for every $x\in E$ and $\phi\in
\mathcal{B}^{+}_{b}(E)$,
\begin{eqnarray}
\mathbb{E}_{x}\left[\mathbb{X}_{t}(\phi
h)\log^{+}\mathbb{X}_{t}(\phi
h)\right]&\le&e^{-\lambda_{1}t}\|\phi\|_{\infty}h(x)\mathbb{E}_{x}\left[e^{\lambda_{1}t}h(x)^{-1}\mathbb{X}_{t}(h)\log^{+}\mathbb{X}_{t}(\phi h)\right]\nonumber\\
&=&e^{-\lambda_{1}t}\|\phi\|_{\infty}h(x)\widetilde{\mathbb{Q}}_{x}\left[\log^{+}\mathbb{X}_{t}(\phi
h)\right].\label{4.2}
\end{eqnarray}
Recall that $\widetilde{\mathcal{G}}_{\infty}$
contains all the information
about the spine.
Since for any $a,b\ge 0$
\begin{equation}
\log^{+}(a+b)\le \log^{+}a+\log^{+}b+\log 2,\label{2.1.1}
\end{equation}
we have by Jensen's inequality
\begin{eqnarray}
\widetilde{\mathbb{Q}}_{x}\left[\log^{+}\mathbb{X}_{t}(\phi
h)\right]&=&\widetilde{\mathbb{Q}}_{x}\left[\widetilde{\mathbb{Q}}_{x}
\left[
\log^{+}\mathbb{X}_{t}(\phi
h) \,\big|\, \widetilde{\mathcal{G}}_{\infty}
\right]
\right]\nonumber\\
&\le&\widetilde{\mathbb{Q}}_{x}\left[\log\widetilde{\mathbb{Q}}_{x}\left[\mathbb{X}_{t}(\phi
h)\vee 1\,\big|\,\widetilde{\mathcal{G}}_{\infty}\right]\right]\nonumber\\
&\le&\widetilde{\mathbb{Q}}_{x}\left[\log^{+}\left(\widetilde{\mathbb{Q}}_{x}\left[\mathbb{X}_{t}(\phi
h)\,\big|\,\widetilde{\mathcal{G}}_{\infty}\right]+1\right)\right]\nonumber\\
&\le&\widetilde{\mathbb{Q}}_{x}\left[\log^{+}\widetilde{\mathbb{Q}}_{x}\left[\mathbb{X}_{t}(\phi
h)\,\big|\,\widetilde{\mathcal{G}}_{\infty}\right]\right]+\log 2.\label{4.4}
\end{eqnarray}
By \eqref{4.2} and \eqref{4.4},
it suffices
to show that
\begin{equation}
\int_{E}\widetilde{\mathbb{Q}}_{x}\left[\log^{+}\widetilde{\mathbb{Q}}_{x}
\left[
\mathbb{X}_{t}(\phi
h)\,\big|\,\widetilde{\mathcal{G}}_{\infty}
\right]
\right]h(x)^{2}m(dx)<+\infty.\label{3.5}
\end{equation}
We get the spine decomposition for
$\widetilde{\mathbb{Q}}_{x}\left[\mathbb{X}_{t}(\phi
h)\,\big|\,\widetilde{\mathcal{G}}_{\infty}\right]$ as follows:
\begin{eqnarray}
\widetilde{\mathbb{Q}}_{x}\left[\mathbb{X}_{t}(\phi
h)\,\big|\,\widetilde{\mathcal{G}}_{\infty}\right]&=&(\phi
h)(\widetilde X_{t})+\sum_{u\prec\xi_{t}}(A_{u}-1)\mathbb{E}_{\widetilde X_{\zeta_{u}}}\left[\mathbb{X}_{t-\zeta_{u}}(\phi
h)\right]\nonumber\\
&\le&(\phi
h)(\widetilde X_{t})+\sum_{u\prec\xi_{t}}A_{u}\|\phi\|_{\infty}e^{-\lambda_{1}(t-\zeta_{u})}h(\widetilde X_{\zeta_{u}}).\label{3.6}
\end{eqnarray}
Note that $\log^{+}(ab)\le \log^{+}a+\log^{+}b$ for any $a,b> 0$.
Using this and an analogy of \eqref{2.1.1} as well as the assumption that $\lambda_1<0$, we have
\begin{eqnarray}
&&\log^{+}\widetilde{\mathbb{Q}}_{x}\left[\mathbb{X}_{t}(\phi
h)\,\big|\,\widetilde{\mathcal{G}}_{\infty}\right]\nonumber\\
&\le&\log^{+}\left(\phi
h(\widetilde X_{t})+\sum_{u\prec\xi_{t}}A_{u}\|\phi\|_{\infty}e^{-\lambda_{1}(t-\zeta_{u})}h(\widetilde X_{\zeta_{u}})\right)\nonumber\\
&\le&\log^{+}(\phi
h(\widetilde X_{t}))+\sum_{u\prec\xi_{t}}\log^{+}\left(A_{u}\|\phi\|_{\infty}e^{-\lambda_{1}(t-\zeta_{u})}h(\widetilde X_{\zeta_{u}})\right)
+ \log^+  n_t \nonumber \\
&\le&\log^{+}(\phi
h(\widetilde X_{t}))+\sum_{u\prec \xi_{t}}\left(
\log^{+}\|\phi\|_{\infty}-\lambda_{1}(t-\zeta_{u})+\log^{+}\left(A_{u}h(\widetilde X_{\zeta_{u}})\right)\right)
+n_{t}  \nonumber\\
&\le&\log^{+}(\phi
h(\widetilde X_{t}))+\sum_{u\prec
\xi_{t}}\log^{+}\left(A_{u}h(\widetilde X_{\zeta_{u}})\right)+
(1 +\log^{+}\|\phi\|_{\infty}-\lambda_{1}t)n_{t}.
\end{eqnarray}
By \eqref{llogl} and using the fact that $\widetilde{m}(dy)$ is an invariant distribution,
$$
\int_{E}\widetilde{\mathbb{Q}}_{x}\left[\log^{+}(\phi
h(\widetilde X_{t}))\right]h(x)^{2}m(dx)=\int_{E}\log^{+}(\phi h)(x)\widetilde{m}(dx)<+\infty.
$$
Hence \eqref{3.5} is implied by
\begin{equation}
\int_{E}\widetilde{\mathbb{Q}}_{x}\left[\sum_{u\prec
\xi_{t}}\log^{+}(A_{u}h(\widetilde X_{\zeta_{u}}))+n_{t}\right]h(x)^{2}m(dx)<+\infty.
\end{equation}
Recall that conditioned on $\mathcal{G}_{\infty}$, $N(\cdot)=\sharp\{i\ge 1:\
(\zeta_{i},A_{i})\in \cdot\}$ is a Poisson random measure on
$[0,+\infty)\times \mathbb{Z}_{+}$ with intensity
$dA^{\widetilde{\mu}}_{s}\sum_{k\in\mathbb{Z}_{+}}\widehat{p}_{k}(\widetilde X_{s})\delta_{k}(dy)$.
We have
\begin{eqnarray}
&&\int_{E}\widetilde{\mathbb{Q}}_{x}\left[\sum_{u\prec
\xi_{t}}\log^{+}(A_{u}h(\widetilde X_{\zeta_{u}}))+n_{t}\right]h(x)^{2}m(dx)\nonumber\\
&=&\int_{E}\widetilde{\mathbb{Q}}_{x}\left[ \widetilde{\mathbb{Q}}_{x}
\left[\sum_{u\prec
\xi_{t}}\log^{+}(A_{u}h(\xi_{\zeta_{u}}))+n_{t}\, \,\Big|\,\mathcal{G}_{\infty}\right]\right]h(x)^{2}m(dx)\nonumber\\
&=&\int_{E}\Pi^{h}_{x}\left[\int_{0}^{t}\sum_{k=0}^{+\infty}\widehat{p}_{k}(\widetilde X_{s})\left(\log^{+}(k
h(\widetilde X_{s}))+1\right)dA^{\widetilde{\mu}}_{s}\right]\widetilde{m}(dx)\nonumber\\
&=&\int_{0}^{t}ds\int_{E}\sum_{k=0}^{+\infty}\widehat{p}_{k}(y)\left(\log^{+}(k
h(y))+1\right)Q(y)h(y)^{2}\mu(dy)\nonumber\\
&\le&t\left(\int_{E}\sum_{k=0}^{+\infty}kp_{k}(y)\log^{+}(k
h(y))h(y)^{2}\mu(dy)+\|Q\|_{\infty}\int_{E}h(y)^{2}\mu(dy)\right).\nonumber
\end{eqnarray}
Immediately the last term is finite by \eqref{llogl}. Hence we complete the proof. \qed

\medskip

\begin{lemma}\label{lem3.2}
If Assumptions 1-2 and \eqref{llogl} hold,
then for any $s,\sigma>0$, $m\in\mathbb{N}$,
and any $x\in E$,
\begin{equation}
\lim_{t\to+\infty}e^{\lambda_{1}(s+t)}\mathbb{X}_{s+t}(\phi
h)-\mathbb{E}_{x}\left[e^{\lambda_{1}(s+t)}\mathbb{X}_{s+t}(\phi
h) \,\big|\, \mathcal{F}_{t}\right]=0\quad\mbox{ in
}L^{1}(\mathbb{P}_{x}),\label{4.9}
\end{equation}
\begin{equation}
\lim_{\mathbb{N}\ni n\to+\infty}e^{\lambda_{1}(n+m)\sigma}\mathbb{X}_{(n+m)\sigma}(\phi
h)-\mathbb{E}_{x}\left[e^{\lambda_{1}(n+m)\sigma}\mathbb{X}_{(n+m)\sigma}(\phi
h) \,\big|\, \mathcal{F}_{n\sigma}\right]=0\quad\mathbb{P}_{x}\mbox{-a.s.}\label{4.10}
\end{equation}
for every $\phi\in \mathcal{B}^{+}_{b}(E)$.
\end{lemma}

\proof For any particle $u\in{\cal Z}_{s}$, let
$\{\mathbb{X}^{u,s}_{t}:t\ge s\}$ denote the branching Markov
process initiated by $u$ at time $s$. It is known that conditioned
on $\mathcal{F}_{s}$, $\mathbb{X}^{u,s}$ and $\mathbb{X}^{v,s}$ are
independent for every $u,v\in{\cal Z}_{s}$ with $u\not=v$. For every
$s,t\ge 0$, define
$$S_{s,t}:=e^{\lambda_{1}t}\mathbb{X}_{s+t}(\phi h)=e^{\lambda_{1}t}\sum_{u\in{\cal Z}_{t}}\mathbb{X}^{u,t}_{s+t}(\phi h),$$
and
$$\widetilde{S}_{s,t}:=e^{\lambda_{1}t}\sum_{u\in{\cal Z}_{t}}\mathbb{X}^{u,t}_{s+t}(\phi h)1_{\{\mathbb{X}^{u,t}_{s+t}(\phi h)\le e^{-\lambda_{1}t}\}}.$$
Obviously $S_{s,t}\ge \widetilde{S}_{s,t}$.

First by the conditional
independence and the Markov property we have
\begin{eqnarray}
&&\mathbb{E}_{x}\left[\left(\widetilde{S}_{m\sigma,n\sigma}-\mathbb{E}_{x}[\widetilde{S}_{m\sigma,n\sigma}
 \,\big|\, \mathcal{F}_{n\sigma}]\right)^{2}\right]\nonumber\\
&=&\mathbb{E}_{x}\left[\mbox{Var}\left[\widetilde{S}_{m\sigma,n\sigma} \,\big|\, \mathcal{F}_{n\sigma}\right]\right]\nonumber\\
&=&e^{2\lambda_{1}n\sigma}\mathbb{E}_{x}\left[\sum_{u\in{\cal Z}_{n\sigma}}\mbox{Var}
\left[\left.\mathbb{X}^{u,n\sigma}_{(n+m)\sigma}(\phi
h)1_{\{\mathbb{X}^{u,n\sigma}_{(n+m)\sigma}(\phi h)\le
e^{-\lambda_{1}n\sigma}\}}\,\right|\,\mathcal{F}_{n\sigma}\right]\right]\nonumber\\
&\le&e^{2\lambda_{1}n\sigma}\mathbb{E}_{x}\left[\left.\sum_{u\in{\cal Z}_{n\sigma}}\mathbb{E}_{x}
\left[\left(\mathbb{X}^{u,n\sigma}_{(n+m)\sigma}(\phi
h)\right)^{2}1_{\{\mathbb{X}^{u,n\sigma}_{(n+m)\sigma}(\phi h)\le
e^{-\lambda_{1}n\sigma}\}}\,\right|\,\mathcal{F}_{n\sigma}\right]\right]\nonumber\\
&=&e^{2\lambda_{1}n\sigma}\mathbb{E}_{x}\left[\sum_{u\in{\cal Z}_{n\sigma}}
\mathbb{E}_{X_{u}(n\sigma)}\left[\left(\mathbb{X}_{m\sigma}(\phi
h)\right)^{2}1_{\{\mathbb{X}_{m\sigma}(\phi h)\le
e^{-\lambda_{1}n\sigma}\}}\right]\right].\label{4.12}
\end{eqnarray}
Let $g_{m\sigma,n\sigma}(y):=\mathbb{E}_{y}\left[\left(\mathbb{X}_{m\sigma}(\phi
h)\right)^{2}1_{\{\mathbb{X}_{m\sigma}(\phi h)\le
e^{-\lambda_{1}n\sigma}\}}\right]$. Immediately we have
$$g_{m\sigma,n\sigma}(y)\le
e^{-\lambda_{1}n\sigma}\mathbb{E}_{y}[\mathbb{X}_{m\sigma}(\phi
h)]\le e^{-\lambda_{1}(n+m)\sigma}\|\phi\|_{\infty}h(y).$$ Thus
$g_{m\sigma,n\sigma}\in L^{2}(E,m)$, and
\begin{equation}
\|g_{m\sigma,n\sigma}\|_{L^{2}(E,m)}\le
e^{-\lambda_{1}(n+m)\sigma}\|\phi\|_{\infty}.\label{3.12}
\end{equation}
Using \eqref{p3}, we continue the estimates in \eqref{4.12} to
get that for $n\in\mathbb{N}$ with $n\sigma>1$,
\begin{eqnarray}
&&\mathbb{E}_{x}\left[\left(\widetilde{S}_{m\sigma,n\sigma}
-\mathbb{E}_{x} \left[\widetilde{S}_{m\sigma,n\sigma} \,\big|\, \mathcal{F}_{n\sigma} \right] \right)^{2}\right]\nonumber\\
&\le&e^{2\lambda_{1}n\sigma}\mathbb{E}_{x}\left[ \mathbb{X}_{n\sigma}(g_{m\sigma,n\sigma})\right] \nonumber\\
&=&e^{\lambda_{1}n\sigma}h(x)P^{h}_{n\sigma}
( g_{m\sigma,n\sigma}/h)(x)  \nonumber\\
&\le&e^{\lambda_{1}n\sigma}h(x)\left(\langle
g_{m\sigma,n\sigma},h\rangle+
e^{\lambda_h /2}
\widetilde{a}_{1}(x)^{1/2}e^{-\lambda_h n\sigma}\|g_{m\sigma,n\sigma}\|_{L^{2}(E,m)}\right)\nonumber\\
&\le&e^{\lambda_{1}n\sigma}h(x)\langle
g_{m\sigma,n\sigma},h\rangle+
e^{\lambda_h /2}
h(x)\widetilde{a}_{1}(x)^{1/2}e^{-\lambda_h n\sigma-\lambda_{1}m\sigma}\|\phi\|_{\infty}.\label{4.14}
\end{eqnarray}
Note that
\begin{eqnarray}
\sum_{n=1}^{+\infty}e^{\lambda_{1}n\sigma}\langle
g_{m\sigma,n\sigma},h\rangle&=&\sum_{n=1}^{+\infty}e^{\lambda_{1}n\sigma}\int_{E}\mathbb{E}_{y}
\left[ \left(\mathbb{X}_{m\sigma}(\phi
h)\right)^{2}1_{\{\mathbb{X}_{m\sigma}(\phi h)\le
e^{-\lambda_{1}n\sigma}\}}\right] h(y)m(dy)\nonumber\\
&\le&\int_{E}h(y)m(dy)\int_{1}^{+\infty}e^{\lambda_{1}(s-1)\sigma}\mathbb{E}_{y}
\left[\left(\mathbb{X}_{m\sigma}(\phi
h)\right)^{2}1_{\{\mathbb{X}_{m\sigma}(\phi h)\le
e^{-\lambda_{1}s\sigma}\}}\right] ds\nonumber\\
&=&\int_{E}h(y)m(dy)\int_{1}^{+\infty}\frac{1}{-\lambda_{1}\sigma
x^{2}}\mathbb{E}_{y}\left[ \left(\mathbb{X}_{m\sigma}(\phi
h)\right)^{2}1_{\{\mathbb{X}_{m\sigma}(\phi h)\le
xe^{-\lambda_{1}\sigma}\}}\right ] dx\nonumber\\
&=&\int_{E}h(y)m(dy)\int_{0}^{+\infty}z^{2}\mathbb{P}_{y}\left(\mathbb{X}_{m\sigma}(\phi
h)\in dz\right)\int_{1\vee
ze^{\lambda_{1}\sigma}}^{+\infty}\frac{1}{-\lambda_{1}\sigma
x^{2}}dx\nonumber\\
&\le&\frac{1}{-\lambda_{1}\sigma}e^{-\lambda_{1}\sigma}\int_{E}h(y)m(dy)\int_{0}^{+\infty}z\mathbb{P}_{y}\left(\mathbb{X}_{m\sigma}(\phi
h)\in dz\right)\nonumber\\
&=&\frac{1}{-\lambda_{1}\sigma}e^{-\lambda_{1}\sigma}\int_{E}h(y)\mathbb{E}_{y}\left[ \mathbb{X}_{m\sigma}(\phi
h)\right] m(dy)\nonumber\\
&\le&\frac{1}{-\lambda_{1}\sigma}e^{-\lambda_{1}\sigma-\lambda_{1}m\sigma}\|\phi\|_{\infty}\int_{E}h(y)^{2}m(dy)<+\infty,\label{4.16}
\end{eqnarray}
where in the second equality above we used
the change of variables $x=e^{-\lambda_1(s-1)\sigma}$.
It follows from \eqref{4.14} and \eqref{4.16} that
$\sum_{n=1}^{+\infty}\mathbb{E}_{x}\left[\left(\widetilde{S}_{m\sigma,n\sigma}-\mathbb{E}_{x} [\widetilde{S}_{m\sigma,n\sigma} \,\big|\, \mathcal{F}_{n\sigma}] \right)^{2}\right]<+\infty$.
Thus by the Borel-Cantelli lemma,
\begin{equation}
\lim_{n\to+\infty}\widetilde{S}_{m\sigma,n\sigma}-\mathbb{E}_{x}\left[ \widetilde{S}_{m\sigma,n\sigma} \,\big|\, \mathcal{F}_{n\sigma}\right] =0\quad
\mathbb{P}_{x}\mbox{-a.s.}\label{4.17}
\end{equation}
Note that for every $n,m\in\mathbb{N}$, we have
\begin{eqnarray}
\mathbb{E}_{x}\left[ S_{m\sigma,n\sigma}-\widetilde{S}_{m\sigma,n\sigma}\right]
&=&\mathbb{E}_{x}\left[\mathbb{E}_{x}[S_{m\sigma,n\sigma} \,\big|\, \mathcal{F}_{n\sigma}]
-\mathbb{E}_{x}[\widetilde{S}_{m\sigma,n\sigma} \,\big|\, \mathcal{F}_{n\sigma}]\right]\nonumber\\
&=&e^{\lambda_{1}n\sigma}\mathbb{E}_{x}\left[\sum_{u\in
{\cal Z}_{n\sigma}}\mathbb{E}_{X_{u}(n\sigma)}\left[ \mathbb{X}_{m\sigma}(\phi
h)1_{\{\mathbb{X}_{m\sigma}(\phi
h)>e^{-\lambda_{1}n\sigma}\}}\right] \right].\label{4.18}
\end{eqnarray}
Let $f_{m\sigma,n\sigma}(y):=\mathbb{E}_{y}\left[\mathbb{X}_{m\sigma}(\phi
h)1_{\{\mathbb{X}_{m\sigma}(\phi h)>e^{-\lambda_{1}n\sigma}\}}\right]$. Obviously,
\begin{equation}\label{new-eq}
f_{m\sigma,n\sigma}(y)\le \|\phi\|_{\infty}\mathbb{E}_{y}[\mathbb{X}_{m\sigma}(h)]
\le e^{-\lambda_{1}m\sigma}\|\phi\|_{\infty}h(y)
\end{equation}
and so
$f_{m\sigma,n\sigma}\in
L^{2}(E,m)$. Then by \eqref{p3}, we have for $n\in\mathbb{N}$ with $n\sigma>1$,
\begin{eqnarray}
\mbox{RHS of }
\eqref{4.18}&=&e^{\lambda_{1}n\sigma}\mathbb{E}_{x}[\mathbb{X}_{n\sigma}(f_{m\sigma,n\sigma})]
=h(x)P^{h}_{n\sigma}
(f_{m\sigma,n\sigma}/h)(x)\nonumber\\
&\le&h(x)\langle
f_{m\sigma,n\sigma},h\rangle+
e^{\lambda_h /2}
\widetilde{a}_{1}(x)^{1/2}h(x)e^{-\lambda_h n\sigma}\|f_{m\sigma,n\sigma}\|_{L^{2}(E,m)}\nonumber\\
&\le&h(x)\langle
f_{m\sigma,n\sigma},h\rangle+
e^{\lambda_h /2}
\widetilde{a}_{1}(x)^{1/2}h(x)e^{-\lambda_h n\sigma-\lambda_{1}m\sigma}\|\phi\|_{\infty},\label{4.19}
\end{eqnarray}
where in the last inequality we used \eqref{new-eq}. Note that
\begin{eqnarray}
\sum_{n=1}^{+\infty}\langle
f_{m\sigma,n\sigma},h\rangle&=&\sum_{n=1}^{+\infty}\int_{E}\mathbb{E}_{y}\left[\mathbb{X}_{m\sigma}(\phi
h)1_{\{\mathbb{X}_{m\sigma}(\phi
h)>e^{-\lambda_{1}n\sigma}\}}\right]h(y)m(dy)\nonumber\\
&\le&\int_{0}^{+\infty}ds\int_{E}\mathbb{E}_{y}\left[\mathbb{X}_{m\sigma}(\phi
h)1_{\{\mathbb{X}_{m\sigma}(\phi
h)>e^{-\lambda_{1}\sigma s}\}}\right]h(y)m(dy)\nonumber\\
&=&\int_{E}h(y)m(dy)\int_{0}^{+\infty}ds\int_{e^{-\lambda_{1}\sigma
s}}^{+\infty}z\mathbb{P}_{y}\left(\mathbb{X}_{m\sigma}(\phi h)\in
dz\right)\nonumber\\
&=&\int_{E}h(y)m(dy)\int_{1}^{+\infty}\frac{1}{-\lambda_{1}\sigma
t}dt\int_{t}^{+\infty}z\mathbb{P}_{y}\left(\mathbb{X}_{m\sigma}(\phi
h)\in
dz\right)\nonumber\\
&=&\frac{1}{-\lambda_{1}\sigma}\int_{E}h(y)m(dy)\int_{1}^{+\infty}z\mathbb{P}_{y}\left(\mathbb{X}_{m\sigma}(\phi
h)\in
dz\right)\int_{1}^{z}\frac{1}{t}dt\nonumber\\
&=&\frac{1}{-\lambda_{1}\sigma}\int_{E}h(y)m(dy)\int_{1}^{+\infty}z\log^{+}z\mathbb{P}_{y}\left(\mathbb{X}_{m\sigma}(\phi
h)\in dz\right)\nonumber\\
&=&\frac{1}{-\lambda_{1}\sigma}\int_{E}\mathbb{E}_{y}\left[\mathbb{X}_{m\sigma}(\phi
h)\log^{+}\mathbb{X}_{m\sigma}(\phi h)\right]h(y)m(dy).\label{4.30}
\end{eqnarray}
The last term is finite by Lemma \ref{lem2.1}. Thus we get
$$\sum_{n=1}^{+\infty}\mathbb{E}_{x}[S_{m\sigma,n\sigma}-\widetilde{S}_{m\sigma,n\sigma}]
=\sum_{n=1}^{+\infty}\mathbb{E}_{x}\left[\mathbb{E}_{x}[S_{m\sigma,n\sigma} \,\big|\, \mathcal{F}_{n\sigma}]-\mathbb{E}_{x}[\widetilde{S}_{m\sigma,n\sigma} \,\big|\, \mathcal{F}_{n\sigma}]\right]<+\infty.$$
Recall that $S_{m\sigma,n\sigma}\ge \widetilde{S}_{m\sigma,n\sigma}$ for
every $n\ge 0$. Again by the Borel-Cantelli lemma, we have
\begin{equation}
\lim_{n\to
+\infty}S_{m\sigma,n\sigma}-\widetilde{S}_{m\sigma,n\sigma}
=\lim_{n\to+\infty}\mathbb{E}_{x}[S_{m\sigma,n\sigma} \,\big|\, \mathcal{F}_{n\sigma}]-\mathbb{E}_{x}[\widetilde{S}_{m\sigma,n\sigma} \,\big|\, \mathcal{F}_{n\sigma}]=0\quad
\mathbb{P}_{x}\mbox{-a.s.}\label{4.21}
\end{equation}
Now \eqref{4.10} follows from \eqref{4.17} and \eqref{4.21}.

Substituting $n\sigma$ by $t$ and $m\sigma$
by $s$ in \eqref{4.14}, we get for any $t>1$,
$$\mathbb{E}_{x}\left[\left(\widetilde{S}_{s,t}-\mathbb{E}_{x}\left[\widetilde{S}_{s,t} \,\big|\, \mathcal{F}_{t}\right]\right)^{2}\right]
\le e^{\lambda_{1}t}h(x)\langle g_{s,t},h\rangle+
e^{\lambda_h /2}
h(x)\widetilde{a}_{1}(x)^{1/2}e^{-\lambda_h t-\lambda_{1}s}\|\phi\|_{\infty}.$$
Using a similar calculation as in \eqref{4.16}, we get for $s>0$, $\int_{1}^{+\infty}e^{\lambda_{1}t}\langle g_{s,t},h\rangle dt<+\infty$. Consequently $\lim_{t\to+\infty} e^{\lambda_{1}t}\langle g_{s,t},h\rangle=0$. Hence $\mathbb{E}_{x}\left[\left(\widetilde{S}_{s,t}-\mathbb{E}_{x}\left[\widetilde{S}_{s,t} \,\big|\, \mathcal{F}_{t}\right]\right)^{2}\right]\to 0$ as $t\to+\infty$, or equivalently
\begin{equation}
\lim_{t\to+\infty} \left( \widetilde{S}_{s,t}-\mathbb{E}_{x}[\widetilde{S}_{s,t} \,\big|\, \mathcal{F}_{t}] \right) =0\quad\mbox{in }L^{2}(\mathbb{P}_{x}).\label{L1}
\end{equation}
Substituting $n\sigma$ by $t$ and $m\sigma$ by $s$ in \eqref{4.19},
we get for any $t>1$,
\begin{eqnarray}
\mathbb{E}_{x}\left[ S_{s,t}-\widetilde{S}_{s,t}\right]&=&\mathbb{E}_{x}\left[\mathbb{E}_{x}
\left[S_{s,t} \,\big|\, \mathcal{F}_{t}\right]-\mathbb{E}_{x}\left[\widetilde{S}_{s,t} \,\big|\, \mathcal{F}_{t}\right]\right]\nonumber\\
&\le&h(x)\langle
f_{s,t},h\rangle+
e^{\lambda_h /2}
\widetilde{a}_{1}(x)^{1/2}h(x)e^{-\lambda_h t-\lambda_{1}s}\|\phi\|_{\infty}.\nonumber
\end{eqnarray}
By similar calculation as in \eqref{4.30}, we get $\int_{0}^{+\infty}\langle
f_{s,t},h\rangle dt<+\infty$ for all $s>0$. Hence $\lim_{t\to+\infty}\langle
f_{s,t},h\rangle=0$. Thus we have
\begin{equation}
\mathbb{E}_{x}\left[ S_{s,t}-\widetilde{S}_{s,t}\right]=\mathbb{E}_{x}\left[\mathbb{E}_{x}
\left[S_{s,t} \,\big|\, \mathcal{F}_{t}\right]-\mathbb{E}_{x}\left[\widetilde{S}_{s,t} \,\big|\, \mathcal{F}_{t}\right]\right]\to 0,\mbox{ as
}t\to+\infty.\label{L2}
\end{equation}
Therefore \eqref{4.9} follows from \eqref{L1} and \eqref{L2}.\qed

\begin{lemma}\label{lem4.3} Under the conditions of Theorem \ref{them2}, we have
\begin{equation}
\lim_{s\to+\infty}\lim_{t\to+\infty}e^{\lambda_{1}(s+t)}\mathbb{E}_{x}\left[ \mathbb{X}_{s+t}(\phi
h)\,\big|\,\mathcal{F}_{t}\right] =M_{\infty}\langle \phi
h,h\rangle\quad\mbox{ in }L^{1}(\mathbb{P}_{x}),
\end{equation}
for every $\phi\in\mathcal{B}^{+}_{b}(E)$ and every $x\in E$.
\end{lemma}
\proof Recall that $M_{t}$ converges to $M_{\infty}$ in
$L^{1}(\mathbb{P}_{x})$ by Theorem \ref{them1}. It suffices to prove
that
\begin{equation}
\lim_{s\to+\infty}\lim_{t\to+\infty}\mathbb{E}_{x}\left[\left|e^{\lambda_{1}(s+t)}\mathbb{E}_{x}[\mathbb{X}_{s+t}(\phi
h)\,\big|\,\mathcal{F}_{t}]-M_{t}\langle\phi h,h\rangle\right|\right]=0.
\end{equation}
Note that by the Markov property,
\begin{eqnarray}
e^{\lambda_{1}(s+t)}\mathbb{E}_{x}\left[\mathbb{X}_{s+t}(\phi
h)\,\big|\,\mathcal{F}_{t}\right]&=&e^{\lambda_{1}(s+t)}\sum_{u\in{\cal Z}_{t}}\mathbb{E}_{X_{u}(t)}\left[\mathbb{X}_{s}(\phi
h)\right]\nonumber\\
&=&e^{\lambda_{1}t}\mathbb{X}_{t}\left(hP^{h}_{s}(\phi)\right).
\end{eqnarray}
Thus we have for any $s,t>t_{0}$,
\begin{eqnarray}
&&\mathbb{E}_{x}\left[\left|e^{\lambda_{1}(s+t)}\mathbb{E}_{x}[\mathbb{X}_{s+t}(\phi
h)\,\big|\,\mathcal{F}_{t}]-M_{t}\langle\phi
h,h\rangle\right|\right]\nonumber\\
&\le&e^{\lambda_{1}t}\mathbb{E}_{x}\left[\left|\mathbb{X}_{t}\left(hP^{h}_{s}\phi-h\langle
\phi h,h\rangle\right)\right|\right]\nonumber\\
&\le&e^{\lambda_{1}t}\mathbb{E}_{x}\left[\mathbb{X}_{t}\left(\left|hP^{h}_{s}\phi-h\langle
\phi h,h\rangle\right|\right)\right]\nonumber\\
&=&h(x)P^{h}_{t}\left(\left|P^{h}_{s}\phi-\langle \phi
h,h\rangle\right|\right)(x).\nonumber
\end{eqnarray}
Using  \eqref{p3} we continue the estimation above to get:
\begin{eqnarray}
&&\mathbb{E}_{x}\left[\left|e^{\lambda_{1}(s+t)}\mathbb{E}_{x}[\mathbb{X}_{s+t}(\phi
h)\,\big|\,\mathcal{F}_{t}]-M_{t}\langle\phi
h,h\rangle\right|\right]\nonumber\\
&\le&
e^{-\lambda_h (s-t_{0}/2)}
\|\phi\|_{\infty}h(x)P^{h}_{t}(\widetilde{a}_{t_{0}}^{1/2})(x)\nonumber\\
&\le&
e^{-\lambda_h (s-t_{0}/2)}
\|\phi\|_{\infty}h(x)\left[\int_{E}\widetilde{a}_{t_{0}}(y)^{1/2}\widetilde{m}(dy)+
\widetilde{a}_{t_{0}}(x)^{1/2}e^{-\lambda_h (t-t_{0}/2)}
\left(\int_{E}\widetilde{a}_{t_{0}}(y)\widetilde{m}(dy)\right)^{1/2}\right]\nonumber\\
&\to&0\quad\mbox{ as }t\to+\infty,\mbox{ and then } s\to+\infty.\nonumber
\end{eqnarray}\qed

\noindent\textit{Proof of Theorem \ref{them2}:} Let $\phi=f/h$, then $\phi\in\mathcal{B}^{+}_{b}(E)$. Note that for any $s,t>0$
\begin{eqnarray}
\left|e^{\lambda_{1}(t+s)}\mathbb{X}_{t+s}(\phi h)-M_{\infty}\langle \phi h,h\rangle\right|
&\le&\left|e^{\lambda_{1}(t+s)}\mathbb{X}_{t+s}(\phi h)-e^{\lambda_{1}(s+t)}\mathbb{E}_{x}\left[\mathbb{X}_{s+t}(\phi
h)\,\big|\,\mathcal{F}_{t}\right]\right|\nonumber\\
&&+\left|e^{\lambda_{1}(s+t)}\mathbb{E}_{x}\left[\mathbb{X}_{s+t}(\phi
h)\,\big|\,\mathcal{F}_{t}\right]-M_{\infty}\langle \phi h,h\rangle\right|.\nonumber
\end{eqnarray}
Therefore Theorem \ref{them2} follows immediately from Lemma \ref{lem3.2} and Lemma \ref{lem4.3}.\qed

\section{Strong law of large numbers}

In this section, we prove  Theorem \ref{them3}.

\subsection{SLLN along lattice times}

\begin{lemma}\label{lem4}
Suppose the assumptions
of
Theorem \ref{them3} hold.
Then for any $\sigma>0$ and any $x\in E$,
\begin{equation}
\lim_{n\to+\infty}e^{\lambda_{1}n\sigma}\mathbb{X}_{n\sigma}(\phi
h)=M_{\infty}\langle \phi h,h\rangle\quad
\mathbb{P}_{x}\mbox{-a.s.}\label{lattice}
\end{equation}
for every $\phi\in \mathcal{B}^{+}_{b}(E)$
\end{lemma}

\proof Recall that
$M_{\infty}=\lim_{n\to+\infty}e^{\lambda_{1}n\sigma}\mathbb{X}_{n\sigma}(h)$.
Let $g(y):=(\phi h)(y)-h(y)\langle \phi h,h\rangle$. The convergence
in \eqref{lattice} is equivalent to
\begin{equation}
\lim_{n\to+\infty}e^{\lambda_{1}n\sigma}\mathbb{X}_{n\sigma}(g)=0\quad\mathbb{P}_{x}\mbox{-a.s.}\label{4.20}
\end{equation}
For an arbitrary $m\in\mathbb{N}$,
\begin{eqnarray}
&&e^{\lambda_{1}(n+m)\sigma}\mathbb{X}_{(n+m)\sigma}(g)\nonumber\\
&=&\left(e^{\lambda_{1}(n+m)\sigma}\mathbb{X}_{(n+m)\sigma}(g)-
\mathbb{E}_{x}\left[\left.e^{\lambda_{1}(n+m)\sigma}\mathbb{X}_{(n+m)\sigma}(g)\,\right|\,\mathcal{F}_{n\sigma}\right]\right)
+\mathbb{E}_{x}\left[\left.e^{\lambda_{1}(n+m)\sigma}\mathbb{X}_{(n+m)\sigma}(g)\,\right|\,\mathcal{F}_{n\sigma}\right]\nonumber\\
&=:&I_{n}+II_{n}.\label{4.25}
\end{eqnarray}
Note that by Lemma \ref{lem3.2}, we have
\begin{equation}
I_{n}=e^{\lambda_{1}(n+m)\sigma}\mathbb{X}_{(n+m)\sigma}(\phi h)-
\mathbb{E}_{x}\left[\left.e^{\lambda_{1}(n+m)\sigma}\mathbb{X}_{(n+m)\sigma}(\phi
h)\,\right|\,\mathcal{F}_{n\sigma}\right]\to 0\quad\mbox{as
}n\uparrow+\infty\quad\mathbb{P}_{x}\mbox{-a.s.}\label{4.26}
\end{equation}
On the other hand, by Markov property, we have
\begin{eqnarray}
II_{n}&=&\sum_{u\in{\cal Z}_{n\sigma}}e^{\lambda_{1}(n+m)\sigma}\mathbb{E}_{X_{u}(n\sigma)}[\mathbb{X}_{m\sigma}(g)]\nonumber\\
&=& \sum_{u\in
{\cal Z}_{n\sigma}}e^{\lambda_{1}(n+m)\sigma}P^{(Q-1)\mu}_{m\sigma}g(X_{u}(n\sigma)) \nonumber\\
&=&
e^{\lambda_{1}(n+m)\sigma}\mathbb{X}_{n\sigma}(P^{(Q-1)\mu}_{m\sigma}g)
.\label{4.22}
\end{eqnarray}
Note that our assumption \eqref{0.2} implies \eqref{0.3}. Then for any fixed $\varepsilon>0$, there exist $m>0$ sufficiently large such that
$$
\sup_{x,y\in E}
|p^{h}(m\sigma,x,y)
-1|
\le \varepsilon.
$$
Then for any $x\in E$,
\begin{eqnarray}
|e^{\lambda_{1}m\sigma}P^{(Q-1)\mu}_{m\sigma}g(x)|&=&|h(x)P^{h}_{m\sigma}
(g/h)(x) |\nonumber\\
&=&h(x)\left|\int_{E}\left(p^{h}(m\sigma,x,y)-1\right)\phi(y)\widetilde{m}(dy)\right|\nonumber\\
&\le&h(x)\int_{E}\left|p^{h}(m\sigma,x,y)-1\right|\phi(y)h(y)^{2}m(dy)\nonumber\\
&\le&\varepsilon h(x)\langle \phi h,h\rangle.\nonumber
\end{eqnarray}
Consequently by \eqref{4.22}, we have
\begin{equation}
|II_{n}|\le \varepsilon\langle \phi h,h\rangle M_{n\sigma}.\label{4.28}
\end{equation}
It follows from \eqref{4.25}, \eqref{4.26} and \eqref{4.28} that
$$\limsup_{n\to+\infty}|e^{\lambda_{1}n\sigma}\mathbb{X}_{n\sigma}(g)|\le\varepsilon\langle
\phi h,h\rangle M_{\infty}\quad\mathbb{P}_{x}\mbox{-a.s.}$$ Hence
we get \eqref{4.20} by letting $\varepsilon\to 0$.\qed

\subsection{Transition from lattice times to continuous time}

In this subsection we extend the convergence along lattice times in Lemma
\ref{lem4} to convergence
along continuous time
and give a sketch of proof of Theorem
\ref{them3}. The main approach in this subsection is similar to that
of \cite[Theorem 3.7]{Chen&Shiozawa} (see also \cite[Theorem
1']{Asmussen}). According to the proof of \cite[Theorem 3.7]{Chen&Shiozawa}, to prove Theorem
\ref{them3}, it suffices to prove the following lemma:

\begin{lemma}\label{limit-U}  Under the conditions of Theorem \ref{them3}, for every open subset $U$ in $
E$ and every $x\in E$,
\begin{equation}
\liminf_{t\to+\infty}e^{\lambda_{1}t}\mathbb{X}_{t}(1_{U}h)\ge
M_{\infty}\int_{U}h(y)^{2}m(dy)\quad\mathbb{P}_{x}\mbox{-a.s.,}\label{5.1}
\end{equation}
\end{lemma}

Note that if \eqref{5.1} is true for any bounded open set $U$ in $E$ with $\overline U\subset E$, then  \eqref{5.1} is true for any open set $U$. In fact, for an arbitrary open set $U$ in $E$, there exists a sequence of bounded open sets  $\{U_n:n\ge 0\}$ such that $\overline U_n\subset E$ and $U_n\uparrow U$.
Hence if $U_{n}$ satisfies \eqref{5.1}, we can deduce that $U$ satisfies \eqref{5.1} by monotone convergence theorem.

In the following we assume that $U$ is an arbitrary bounded open set in $E$ with $\overline U\subset E$.
For any
$\varepsilon,\sigma>0$, $n\in\mathbb{N}$ and $x\in E$, define
$$U^{\varepsilon}(x):=\{y\in U:\ h(y)\ge \frac{1}{1+\varepsilon}\,h(x)\},$$
$$
Y^{\sigma,\varepsilon}_{n,u}:=\frac{1}{1+\varepsilon}\,(h1_U)(X_{u}(n\sigma))
1_{\left\{X_{v}(t)\in U^{\varepsilon}(X_{u}(n\sigma))\mbox{ for all
}v\in {\cal L}^{u,n\sigma}_{t}\mbox{ and }t\in
[n\sigma,(n+1)\sigma]\right\}},
$$
and
$$
S^{\sigma,\varepsilon}_{n}:=\sum_{u\in{\cal Z}_{n\sigma}}e^{\lambda_{1}n\sigma}Y^{\sigma,\varepsilon}_{n,u}.
$$
Here for each $t\ge n\sigma$, ${\cal L}^{u,n\sigma}_{t}$ denotes
the set of particles which are alive at $t$ and are descendants of the particle $u\in\mathcal{L}_{n\sigma}$.

To prove Lemma \ref{limit-U}, we need the following lemma.
\begin{lemma}\label{lem5.3}
 Suppose that $U$ is an arbitrary bounded open set $U$ in $E$ with $\overline U\subset E$.
Under the conditions of Theorem \ref{them3}, we have
\begin{equation}
\lim_{n\to+\infty}S^{\sigma,\varepsilon}_{n}-\mathbb{E}_{x}\left[S^{\sigma,\varepsilon}_{n}\,\big|\,\mathcal{F}_{n\sigma}\right]=0\quad\mathbb{P}_{x}\mbox{-a.s.}\label{5.2}
\end{equation}
for every $x\in E$.
\end{lemma}

\proof
Note that for $n\in\mathbb{N}$ and $\{Y_{i}:\ i=1,\cdots,n\}$
independent real-valued centered random variables,
$$
\mathrm{E}|\sum_{i=1}^{n}Y_{i}|^{2}
=\mbox{Var}\sum_{i=1}^{n}Y_{i}=\sum_{i=1}^{n}\mathrm{E}|Y_{i}|^{2}.
$$
Thus by the conditional independence between
subtrees and the Markov property of a branching Markov process, we
have
\begin{eqnarray*}
\mathbb{E}_{x}\left[\left.\left(S^{\sigma,\varepsilon}_{n}-
\mathbb{E}_{x}\left[S^{\sigma,\varepsilon}_{n}\,\big|\,\mathcal{F}_{n\sigma}\right]\right)^{2}
\right|\mathcal{F}_{n\sigma}\right]
&=&
e^{2\lambda_{1}n\sigma}\sum_{u\in{\cal Z}_{n\sigma}}\mathbb{E}_{x}
\left[\left.\left|Y^{\sigma,\varepsilon}_{n,u}-\mathbb{E}_{x}
\left[Y^{\sigma,\varepsilon}_{n,u}\,\big|\,\mathcal{F}_{n\sigma}\right]\right|^{2}\right|\mathcal{F}_{n\sigma}\right]\nonumber\\
&\le&2
e^{2\lambda_{1}n\sigma}\sum_{u\in{\cal Z}_{n\sigma}}\mathbb{E}_{x}
\left[|Y^{\sigma,\varepsilon}_{n,u}|^{2}\,\big|\,\mathcal{F}_{n\sigma}\right]\nonumber\\
&=&2
e^{2\lambda_{1}n\sigma}\sum_{u\in{\cal Z}_{n\sigma}}\mathbb{E}_{X_{u}(n\sigma)}
\left[|Y^{\sigma,\varepsilon}_{0,\emptyset}|^{2}\right]\nonumber\\
&\le&2
(1+\varepsilon)^{-2}e^{2\lambda_{1}n\sigma}\sum_{u\in{\cal Z}_{n\sigma}}(h^21_U)
(X_{u}(n\sigma)).
\end{eqnarray*}
Thus
\begin{eqnarray*}
\sum_{n=1}^{+\infty}\mathbb{E}_{x}\left[\left(S^{\sigma,\varepsilon}_{n}-\mathbb{E}_{x}
\left[S^{\sigma,\varepsilon}_{n}\,\big|\,\mathcal{F}_{n\sigma}\right]\right)^2\right]
&\le&2
(1+\varepsilon)^{-2}\sum_{n=1}^\infty  e^{\lambda_{1}2n\sigma}\mathbb{E}_{x}\left[\sum_{u\in{\cal Z}_{n\sigma}}(h^21_U)(X_{u}(n\sigma))\right]\\
&\le&2
(1+\varepsilon)^{-2}\sum_{n=1}^\infty e^{\lambda_{1}2n\sigma}P^{(Q-1)\mu}_{n\sigma}(h^{2}1_U)(x)\\
&\le&2
(1+\varepsilon)^{-2}h(x)\sum_{n=1}^\infty e^{\lambda_{1}2n\sigma}P^{h}_{n\sigma}(h1_U)(x).
\end{eqnarray*}
By Proposition \ref{P:1.8}, for any $\eps >0$,
$|p^{h}(n\sigma,x,y)
-1|
\le\eps $ for every $x,y\in E$
when $n$
is sufficiently large.
It follows that
$P^{h}_{n\sigma}(h1_{U})(x)\le(1+\eps )\int_{U}h(y)\widetilde{m}(dy)<+\infty$,
and so
$$
\sum_{n=1}^{+\infty}\mathbb{E}_{x}\left[\left(S^{\sigma,\varepsilon}_{n}-\mathbb{E}_{x}
\left[S^{\sigma,\varepsilon}_{n}\,\big|\,\mathcal{F}_{n\sigma}\right]\right)^2\right]\le ch(x)\sum_{n=1}^\infty e^{\lambda_{1}n\sigma}<\infty,$$
where $c$ is a positive constant.
This together with the Borel-Cantelli lemma yields \eqref{5.2}.
\qed

\bigskip
\noindent\textit{Proof of Lemma \ref{limit-U} :}  If $U$ is an arbitrary bounded open set
in $E$ with $\overline U\subset E$, then
\eqref{5.1} follows from Lemma \ref{lem5.3} in the same way as \cite[Theorem3.9]{Chen&Shiozawa}.
We omit the details here. By the argument after \eqref{5.1}, we know that \eqref{5.1}
holds  for any open subset $U$ of  $E$.

\bigskip
\noindent\textbf{Acknowledgement} The research of Zhen-Qing Chen is partially supported
by NSF grant DMS-1206276 and NNSFC 11128101. The research of Yan-Xia Ren is supported by NNSFC (Grant No.  11271030 and 11128101). The research of Ting Yang is partially supported by NNSF of China (Grant No. 11501029) and Beijing Institute of Technology Research Fund Program for Young Scholars.

\end{document}